\documentclass[10pt,a4paper,twoside]{article}
\usepackage[usenames,dvipsnames]{xcolor}
\usepackage{hyperref}
\definecolor{DarkPurple}{rgb}{0.40,0.0,0.20}
\hypersetup{colorlinks=true,citecolor=NavyBlue,linkcolor=DarkPurple,urlcolor=NavyBlue, breaklinks=true}
\usepackage{color} 

\usepackage{xstring}
\usepackage{dynkin-diagrams}
\usepackage{etoolbox}
\usepackage[flushleft]{threeparttable}
\usepackage{hhline}
\usepackage{comment}
\usepackage{tikz}
\usepackage[page]{appendix}
\usepackage[utf8]{inputenc}
\usepackage{filecontents}
\usepackage{graphicx} \usepackage{mathtools}
\usepackage{etoolbox}
\usepackage{amsmath}
\usepackage{amsthm}
\usepackage{amsfonts}
\usepackage{amssymb}
\usepackage{mathrsfs}
\usepackage{xypic}
\usepackage{xcolor}
\usepackage{titlesec}
\usepackage{titletoc}
\usepackage{tabularx}
\newtheorem{lem}{Lemma}

\newtheorem{cor}{Corollary}
\newtheorem{prop}{Proposition}
\newtheorem{thm}{Theorem}
\newtheorem*{quest*}{Question}
\newtheorem{quest}{Question}
\newtheorem*{thm*}{Theorem}
\newtheorem{dfn}{Definition}
\DeclareMathAlphabet{\mathpzc}{OT1}{pzc}{m}{it}

\begin{document}

\newcommand{\chicyc}{\chi_{\operatorname{cyc}}}
\newcommand{\Brnr}{\operatorname{Br}_{\operatorname{nr}}}
\newcommand{\Bronr}{\operatorname{Br}_{1,\operatorname{nr}}}
\newcommand{\Br}{\operatorname{Br}}
\newcommand{\LT}{\operatorname{LT}}
\newcommand{\Qp}{\mathbb{Q}_p}
\newcommand{\TODO}{{\color{red} TODO}}
\newcommand{\Gm}{\mathbb{G}_m}
\newcommand{\Scin}{\operatorname{Scin}}
\newcommand{\Fil}{\operatorname{Fil}}
\newcommand{\Ind}{\operatorname{Ind}}
\newcommand{\Sym}{\operatorname{Sym}}
\newcommand{\sym}{\operatorname{sym}}
\newcommand{\semis}[1]{#1\operatorname{-ss}}
\newcommand{\Alt}{\operatorname{Alt}}
\newcommand{\DGK}{\operatorname{DG}_K}
\newcommand{\Img}{\operatorname{Im}}
\newcommand{\Ker}{\operatorname{Ker}}
\newcommand{\grp}{\operatorname{grp}}
\newcommand{\cont}{\operatorname{cont}}
\newcommand{\ord}{\operatorname{ord}}
\newcommand{\Hom}{\operatorname{Hom}}
\newcommand{\Ext}{\operatorname{Ext}}
\newcommand{\Ad}{\operatorname{Ad}}
\newcommand{\Id}{\operatorname{id}}
\newcommand{\Lie}{\operatorname{Lie}}
\newcommand{\Lift}{\operatorname{Lift}}
\newcommand{\ad}{\operatorname{ad}}
\newcommand{\rk}{\operatorname{rk}}
\newcommand{\Det}{\operatorname{Det}}
\newcommand{\LHS}{\operatorname{LHS}}
\newcommand{\RHS}{\operatorname{RHS}}
\newcommand{\bFp}{\bar{\mathbb{F}}_p}
\newcommand{\bZp}{\bar{\mathbb{Z}}_p}
\newcommand{\bQp}{\bar{\mathbb{Q}}_p}
\newcommand{\bF}{{\mathbb{F}}}
\newcommand{\bZ}{{\mathbb{Z}}}
\newcommand{\bA}{{\mathbb{A}}}
\newcommand{\bB}{{\mathbb{B}}}
\newcommand{\bQ}{{\mathbb{Q}}}
\newcommand{\bR}{{\mathbb{R}}}
\newcommand{\bC}{{\mathbb{C}}}
\newcommand{\scrC}{{\mathscr{C}}}
\newcommand{\Fq}{\bar{\mathbb{F}}_q}
\newcommand{\GSp}{\operatorname{GSp}}
\newcommand{\Sp}{\operatorname{Sp}}
\newcommand{\GL}{\operatorname{GL}}
\newcommand{\cris}{\operatorname{cris}}
\newcommand{\modp}{\bF}
\newcommand{\End}{\operatorname{End}}
\newcommand{\SL}{\operatorname{SL}}
\newcommand{\SO}{\operatorname{SO}}
\newcommand{\GO}{\operatorname{GO}}
\newcommand{\cO}{\mathcal{O}}
\newcommand{\Lieg}{\mathfrak{g}}
\newcommand{\Lieh}{\mathfrak{h}}
\newcommand{\Mat}{\operatorname{Mat}}
\newcommand{\Gal}{\operatorname{Gal}}
\newcommand{\Art}{\operatorname{Art}}
\newcommand{\BHT}{\mathbb{B}_{\operatorname{HT}}}
\newcommand{\DXT}{{\hatI}^{X(T)}}
\newcommand{\HTGC}{\prod_{\sigma:K_m\hookrightarrow \bC}X_*(G_{\bC})}
\newcommand{\HTGCa}{\prod_{\tilde\tau\in S}X_*(G_{\bC})}
\newcommand{\HT}{{\operatorname{HT}}}
\newcommand{\hT}{{\mathcal{HT}}}
\newcommand{\dR}{{\operatorname{dR}}}
\newcommand{\hatbI}{\widehat{\bar I}_K}
\newcommand{\hatI}{\widehat{I}_K}
\newcommand{\symt}{\operatorname{sym}^2}
\newcommand{\altt}{\operatorname{alt}^2}
\newcommand{\socle}{\operatorname{soc}}
\newcommand\nlcup{\mathrel{\ooalign{\hss$\cup$\hss\cr \kern0.7ex\raise0.6ex\hbox{\scalebox{0.4}{$\diamond$}}}}}
\newcommand{\Dynkin}{\operatorname{Dynkin}}
\newcommand{\Siegal}{\operatorname{Siegal}}
\newcommand{\Aut}{\operatorname{Aut}}
\newcommand{\red}{\operatorname{red}}
\newcommand{\ur}{\operatorname{ur}}
\newcommand{\urind}{\operatorname{urind}}
\newcommand{\Isom}{\operatorname{Isom}}
\newcommand{\IsomFil}{\operatorname{Isom-fil}}
\newcommand{\crys}{\operatorname{crys}}
\newcommand{\can}{\operatorname{can}}
\newcommand{\val}{\operatorname{val}_p}
\newcommand{\irr}{\operatorname{irr}}
\newcommand{\sems}{\operatorname{ss}}
\newcommand{\std}{\operatorname{std}}
\newcommand{\Vect}{\operatorname{Vect}}
\newcommand{\der}{\operatorname{der}}
\newcommand{\Ord}{\operatorname{Ord}}
\newcommand{\gr}{\operatorname{gr}}
\newcommand{\Res}{\operatorname{Res}}
\newcommand{\res}{\operatorname{res}}
\newcommand{\Rep}{\operatorname{Rep}}
\newcommand{\Tor}{\operatorname{Tor}}
\newcommand{\ab}{\operatorname{ab}}
\newcommand{\noyaux}[1]{\operatorname{ker}{#1}}
\newcommand{\spa}{\operatorname{Spa}}
\newcommand{\spec}{\operatorname{Spec}}
\newcommand{\spf}{\operatorname{Spf}}
\newcommand{\Span}{\operatorname{span}}
\newcommand{\opn}[1]{\operatorname{#1}}
\newcommand{\dirlim}[1]{\underset{#1}{\underrightarrow{\operatorname{lim}}}}
\newcommand{\invlim}[1]{\underset{#1}{\underleftarrow{\operatorname{lim}}}}

\title{Crystalline Lifts and a Variant of the Steinberg--Winter Theorem}
\date{}
\author{Zhongyipan Lin}
\maketitle
\begin{abstract}
Let $K/\Qp$ be a finite extension.
For all irreducible representations
$\bar\rho: G_K \to G(\bFp)$
valued in a general reductive group $G$,
we construct crystalline lifts
of $\bar\rho$ which are Hodge--Tate regular.
We also discuss rationality questions.
We prove a variant of the Steinberg--Winter theorem
along the way.
\end{abstract}

\section{Introduction}~

Fix a connected split reductive group $G$.

\subsection{}

It is often desirable to describe
automorphisms of a reductive group
in root-theoretic terms.
When we are concerned with finite order automorphisms
of a reductive group over a characteristic $0$ field,
the following theorem suffices
\begin{thm*} [Steinberg-Winter, {\cite[Theorem 7.5]{St68}}]
Let $M$ be a linear algebraic group over a field $k$.
Let $\bar k$ be the algebraic closure of $k$.
Let $F_M : M\to M$ be an automorphism of $M$ which can be realized
as conjugation by an element $g\in G(k)$ after an embedding
$M \hookrightarrow G$.

If $g$ is a semi-simple element, then $F_M$ fixes a maximal torus of $M_{\bar k}$.
\end{thm*}
However, in prime characteristic, finite order automorphisms are not always semi-simple, and we need a new criterion for the existence of maximal tori fixed by $F_M$.

Note that the semi-simplicity of $F_M$ implies that the subgroup $\Gamma_{F_M}$
generated by $g$ and $Z_M(M)^\circ$ (the neutral component of the center of $M$)
is a $G$-completely reducible subgroup of $G(\bar k)$.
We conjecture that the $G$-complete reducibility 
of $\Gamma_{F_M}$ is sufficient for (and, to a certain degree, characterizes)
the existence of an $F_M$-fixed maximal torus of $M$.

The notion of $G$-complete reducibility was
introduced by Serre.
In his Moursund Lectures \cite{Se98},
a subgroup $\Gamma \subset G(\bar k)$ is defined to be \textit{$G$-completely reducible}
if for any parabolic subgroup $P$ of $G_{\bar k}$ containing $\Gamma$,
    a Levi subgroup of $P$ also contains $\Gamma$. 
Similarly, we say
a subgroup $\Gamma \subset G(\bar k)$ is \textit{$G$-irreducible}
if it is not contained in any proper parabolic subgroup of $G$.

We prove the following:

\begin{thm} [\ref{thm:sw-irr}, \ref{cor:sw-max-rank}]
Let $M$ be a connected reductive group over a field $k$.
Let $\bar k$ be the algebraic closure of $k$.
Let $F_M : M\to M$ be an automorphism of $M$ which can be realized
as conjugation by an element $g\in G(k)$ after an embedding
$M \hookrightarrow G$.

Let $\Gamma_{F_M}$ be the subgroup of $G$ generated by $g$ and $Z_M(M)^\circ$.
If either
\begin{itemize}
\item $\Gamma_{F_M}$ is $G$-irreducible or
\item $\Gamma_{F_M}$ is $G$-completely reducible, $\rk M = \rk G$, $\operatorname{char} k \ne 2~\text{or}~3$, and $M$ has connected center.
\end{itemize}
then $F_M$ fixes a maximal torus of $M_{\bar k}$.
\end{thm}
\noindent
which suffices for our application to the theory of Galois representations.

We believe our new method can be used to establish a stronger form of Steinberg-Winter by
developing the theory of $G$-complete reducibility for (possibly disconnected) linear algebraic groups
using dynamic methods (\ref{dyn-fil}).
We don't pursue this because we don't want to digress too much.
We do explain how this can possibly be done.

\subsection{}
Let $K/\Qp$ be a finite extension.
We are interested in the following question:
\begin{quest}
\label{quest:lift}
Let $\bar \rho:G_K\to G(\bFp)$ be a group homomorphism.
Does there exist a crystalline representation
$\rho:G_K\to G(\bZp)$ such that $\rho \equiv \bar\rho$?
\end{quest}
\noindent
Question \ref{quest:lift} is raised in \cite{CL11} for $G=\GL_N$,
where they used the machinery of $p$-adic Hodge theory 
to study general torsion Galois representations;
it also has global applications such as
constructing geometric Galois representations
(see, for example, \cite{FKP21}) which conjecturally
correspond to algebraic automorphic forms.

Any characteristic $p$ representation of $G_K$
is an extension of $G$-completely reducible representations.
To construct crystalline lifts of general characteristic~$p$
representations of $G_K$,
it is a common strategy (for example \cite{EG19} for $\GL_N$ and \cite{L21} for $G_2$ and classical groups)
to first construct
lifts of $G$-completely reducible representations, and then try to lift the extension class.

In this paper, we carry out the first steps of the above strategy.
We prove the following theorem:
\begin{thm} [\ref{thm:irreducible-lift}, \ref{thm:ht-regular-lift}, \ref{cor:lift-witt}]
Let $\kappa$ be the residue field of $K$. 
Let $K^{\ur}$ be the maximal unramified extension of $K$ in a fixed algebraic closure.
Let $\bF/\kappa$ be a finite extension of degree $f$.

Let $\bar \rho : G_K \to G(\bF)$ be a group homomorphism whose image is 
a $G$-completely reducible subgroup.
Assume $G$ is split.

\begin{itemize}
\item There exists a characteristic $0$ lift $\rho: G_K\to G(W(\bFp))$ of $\bar\rho$;
\item There exists a Hodge-Tate regular crystalline lift $\rho: G_K\to G(K^{\ur})$ of $\bar\rho$.
\end{itemize}
\end{thm}

We discuss Hodge-Tate theory of Galois representations valued in general reductive groups
in section \ref{sec:hodge-tate-cocharacter}.

\subsection*{Acknowledgement}
It is an immense pleasure to express my deep sense of gratitude
to
my PhD advisor David Savitt,
for suggesting to me the question of constructing
    crystalline lifts of Galois representations
    valued in general reductive groups, and
    for reading the draft carefully and making
    numerous helpful comments and suggestions.

We thank the anonymous referee for a very careful reading and numerous corrections and suggestions.

\section{A variant of Steinberg-Winter theorem}

The key tool in this section is dynamic methods.

\subsection{Dynamic methods}
\label{dyn-fil}
We review \cite[Section 4.1]{Crd14}.
Let $X$ be a scheme over a base scheme $S$, and fix a $\Gm$-action 
$m:\Gm\times X \to X$ on $X$.
For each $x\in X(S)$, we say
$$
\lim_{t\to 0} m(t, x) \hspace{3mm}\text{exists,}
$$
if the morphism $\Gm \to X$, $t\mapsto m(t,x)$
extends a a morphism $\bA^1\to X$.
If the limit exists, the origin $0 \in \bA^1(S)$ maps
to a unique element $\alpha \in X(S)$;
we write $\lim_{t\to 0} m(t, x) = \alpha$.

Let $\lambda$ be a cocharacter of a reductive group $G$ over a field $k$.
Define the following functor on the category of $k$-algebras
$$P_G(\lambda)(A)=\{g\in G(A)|\lim_{t\to 0}\lambda(t)g\lambda(t)^{-1}\text{~exists.}\}$$
where $A$ is a general $k$-algebra. 

Define $$U_G(\lambda)(A)=\{g\in G(A)|\lim_{t\to 0}\lambda(t)g\lambda(t)^{-1}=1\},$$
and denote by $Z_G(\lambda)$ the centralizer of $\lambda$ in $G$.

Since $G$ is a reductive group over a field, $P_G(\lambda)$ is a parabolic subgroup of $G$,
$U_G(\lambda)$ is the unipotent radical of $P_G(\lambda)$,
and $Z_G(\lambda)$ is a Levi subgroup of $P_G(\lambda)$.

The following proposition is the first application of dynamic methods in this section,
and motivates us to consider $G$-compete reducibility in Steinberg-Winter
type questions.

\begin{prop}
Let $M$ be a connected reductive group over a field $k$.
Let $\bar k$ be the algebraic closure of $k$.
Let $F_M : M\to M$ be an automorphism of $M$ which can be realized
as conjugation by an element $g\in G(k)$ after an embedding
$M \hookrightarrow G$.

If $g$ is semisimple, then
$g$ and $Z_M(M)^\circ$ generate a $G$-completely reducible subgroup.
\end{prop}

\begin{proof}
Let $P \subset G_{\bar k}$ be a parabolic subgroup of $G$ which
contains both $g$ and $Z_M(M)^\circ$.
We want to show a Levi subgroup of $P$
also contains both $g$ and $Z_M(M)^\circ$.

Put $L := Z_P(Z_M(M)^\circ)$, the centralizer of $Z_M(M)^\circ$
in $P$.
Note that conjugation by $g$ fixes $L$.
We claim $L$ contains a maximal torus of $G$.
Since $Z_M(M)^\circ$ is a torus, 
it is contained in a maximal torus of $P$.
A maximal torus of $P$ is also a maximal torus of $G$.
Any maximal torus containing $Z_M(M)^\circ$ 
is in the centralizer of $Z_M(M)^\circ$ because of commutativity of tori.

By Steinberg-Winter, there exists a maximal torus $T \subset L$
which is fixed by $g$. By the previous paragraph, $T$ is also a maximal torus
of $G$.

By dynamic methods, there exists a cocharacter $\lambda:\Gm \to T$
such that $P = P_G(\lambda)$.
The two cocharacters $\lambda, g\lambda g^{-1}: \Gm \to T$
lie in the same maximal torus, and can be regarded as elements of the
cocharacter lattice $X_*(G, T)$.
Since $g\in P$, $g(P_G(\lambda))g^{-1} = P_G(g\lambda g^{-1})=P_G(\lambda)$.
So $\lambda, g\lambda g^{-1} \in X_*(G,T)$ are in the same Weyl chamber.
Since $g \in N_G(T)$, $g \lambda g^{-1}$ and $\lambda$ are in the same Weyl orbit,
and thus we must have $\lambda = g \lambda g^{-1}$.
So $g \in Z_G(\lambda)$ and $Z_M(M)^\circ \subset T \subset Z_G(\lambda)$.
Since $Z_G(\lambda)$ is a Levi subgroup of $P$, we are done.
\end{proof}

\subsection{A generalization of dynamic methods}
Dynamic methods allow us to prove theorems over general base schemes by doing 
mathematical analysis.
To do so, we need to generalize the functors $P_G(-)$.

Let $f: \Gm \to G$ be a \textit{$k$-scheme morphism}.
Define the following functor on the category of $k$-algebras
$$P_G(f)(A)=\{g\in G(A)|\lim_{t\to 0}f(t)g f(t)^{-1}\text{~exists.}\}$$
where $A$ is a general $k$-algebra. 
We call $f$ a \textit{fake cocharacter}.
Here ``a limit exists'' means
the scheme morphism  $\Gm \to G$, defined by $t\mapsto f(t)gf(t)^{-1}$,
extends to a scheme morphism $\bA^1\to G$.
Note that $P_G(f)$ is not representable in general.
We define similarly $U_G(f)$.

\begin{lem}
\label{lem:more-dynamic}
Let $G$ be a connected reductive group over a field $k$.
Let $\lambda,\mu:\Gm\to G$
be cocharacters of $G$.
Assume $P_G(\lambda) = P_G(\mu)=:B$
is a Borel subgroup of $G$.
Let~$U$ be the unipotent radical of $B$.

(1)
The functor
$P_G(\mu\lambda)$ is representable by a Borel subgroup.
In fact, we have $P_G(\mu\lambda) = P_G(\mu)=P_G(\lambda)$.

(2)
The limit
$$
\lim _{t\to 0} \lambda(t)\mu(t)\lambda(t)^{-1}\mu(t)^{-1}
$$
exists in the sense of subsection \ref{dyn-fil}, and lies in $U$.

(3)
Let $u$ be an element of $U$.
The limit
$$
\lim _{t\to 0} \lambda(t)u \mu(t) u^{-1}\lambda(t)^{-1}\mu(t)^{-1}
$$
exists in the sense of subsection \ref{dyn-fil} and lies in $U$.

(4)
Now assume $\lambda$ is a product of cocharacters $\lambda_1$, \ldots, $\lambda_s$
such that $P_G(\lambda_i)=B$ for all $i$.
Then $P_G(\lambda)=B$, and the limits in (2) and (3)
still exist and lie in $U$.

Moreover, for any embedding $G \hookrightarrow H$ of connected reductive groups,
$P_H(\lambda)$ is representable by a parabolic subgroup of $H$.
\end{lem}

\begin{proof}
(1)
Since all maximal tori in $B$ are conjugate to each other, there exists 
an element $x \in U_G(\lambda) = U_G(\mu)$ such that 
conjugation by $x$ maps the maximal torus containing $\lambda$
to the maximal torus containing $\mu$.
In particular, $(x \lambda x^{-1}) \mu = \mu (x\lambda x^{-1})$.
Write $\xi$ for $x \lambda x^{-1}$.
We have 
\begin{align*}
&\lim_{t\to 0} \mu(t) \lambda(t) g \lambda(t)^{-1}\mu(t)^{-1}  = 
\lim_{t\to 0} \mu(t) x^{-1}\xi(t) x g  x^{-1} \xi(t)^{-1} x \mu(t)^{-1} \\
&= \lim_{t\to 0} (\mu(t) x^{-1} \mu(t)^{-1})\cdot(\mu(t)\xi(t) x g  x^{-1} \xi(t)^{-1} \mu(t)^{-1})\cdot (\mu(t) x \mu(t)^{-1}) \\
&= \lim_{t\to 0} \mu(t)\xi(t) x g  x^{-1} \xi(t)^{-1} \mu(t)^{-1}
\end{align*}
Note that the last step is because $x\in U_G(\mu)$, and $\lim_{t\to 0} \mu(t) x \mu(t)^{-1}=1$.
So we have $P_G(\mu\lambda) = x^{-1}P_G(\mu\xi)x$.
Since $\mu\xi$ is a genuine cocharacter, $P_G(\mu\xi)$ is representable by a parabolic.

Since $\mu\xi =\xi \mu$, we can regard $\mu$ and $\xi$ as elements in a cocharacter lattice $X_*(G,T)$
where $T$ is a maximal torus containing $\mu$ and $\xi$.
Since $P_G(\mu) = P_G(\lambda) = P_G(\xi)$, $\mu$ and $\xi$ lie in the (interior of the) same Weyl chamber.
The cocharacter $\mu\xi$ is the sum of $\mu$ and $\xi$ in the cocharacter lattice $X_*(G,T)$,
and lies in the same Weyl chamber.
So $P_G(\mu\xi) = P_G(\mu) = P_G(\lambda)$.
Since $x \in B$, we have
$P_G(\mu\lambda) = x^{-1}P_G(\mu\xi)x = P_G(\mu)=P_G(\lambda)$.

(2)
Since all maximal tori of $B$ are conjugate to each other,
there exists an element $g\in U$ such that
$g Z_G(\lambda)g^{-1} = Z_G(\mu)$.
Write $\xi:=g \lambda g^{-1}$, and we have
$\xi \mu = \mu \xi$.
By part (1), $P_G(\xi\mu) = B$.
By the dynamic description of the Borel $B$, the limits
$$\lim_{t\to 0}\xi(t)g \xi(t)^{-1} = 1,$$
$$\lim_{t\to 0}\xi(t)\mu(t)g^{-1} \mu(t)^{-1}\xi(t)^{-1} = 1\text{, and}$$
$$\lim_{t\to 0}\mu(t)g \mu(t)^{-1}=1$$ all exist.
The expression
\begin{align*}
&\lambda(t)\mu(t)\lambda(t)^{-1}\mu(t)^{-1}
=
g^{-1}\xi(t) g \mu(t) g^{-1} \xi(t)^{-1} g \mu(t)^{-1} \\
&=
g^{-1}\cdot (\xi(t) g \xi(t)^{-1}) \cdot (\xi(t) \mu(t) g^{-1} \mu(t)^{-1}\xi(t)^{-1})
\cdot (\mu(t) g \mu(t)^{-1})
\end{align*}
has a limit as $t\to 0$.

(3)
We have
\begin{align*}
&\lambda(t)u \mu(t) u^{-1}\lambda(t)^{-1}\mu(t)^{-1}
 \\=&
(\lambda(t)u \mu(t) u^{-1}\lambda(t)^{-1} u\mu(t)^{-1} u^{-1}) (u \mu(t) u^{-1} \mu(t)^{-1}).
\end{align*}
So (3) follows from (2).

(4)
The method is the same but notations are more complicated.
We define inductively cocharacters $\xi_i$that commute with each other, and elements $u_i$ of~$U$.
Our induction assumption is $P_G(\lambda_1 \cdots \lambda_{j})=P_G(\xi_1\cdots\xi_{j})=B$ for all $j < s$.
Define $\xi_1:=\lambda_1$ and $u_1:=1$.
Let $u_i$ be an element of $U$
such that $\xi_i:=u_i\lambda_iu_i^{-1}$
commutes with $\xi_1\cdots\xi_{i-1}$.
Write $\zeta_j$ for $\xi_1 \xi_2 \cdots\xi
_j$, and write $v_j$ for $u_j/u_{j-1}$
(set $u_0=1$).
We have, for $g\in G$,
\begin{align*}
  \lambda(t) g \lambda(t)^{-1} 
= &(\zeta_1(t) v_2 \zeta_1(t)^{-1})
(\zeta_2(t) v_3 \zeta_2(t)^{-1})\cdots \\
& (\zeta_s(t)u_sgu_s^{-1}\zeta_s(t)^{-1})\\
&(\zeta_{s-1}(t) v_s \zeta_{s-1}(t)^{-1})^{-1}\cdots
(\zeta_1(t) v_2 \zeta_1(t)^{-1})^{-1}
\hspace{10mm}
\text{($\dagger$)}
\end{align*}
which has a limit if and only if $g\in B$.
Similarly, 
\begin{align*}
  &\lambda(t) \mu(t) \lambda(t)^{-1}\mu(t)^{-1} 
\\=& (\zeta_1(t) v_2 \zeta_1(t)^{-1})
(\zeta_2(t) v_3 \zeta_2(t)^{-1})\cdots \\
& (\zeta_s(t)u_s \mu(t) u_s^{-1}\zeta_s(t)^{-1} \mu(t)^{-1}\\
&\mu(t)(\zeta_{s-1}(t) v_s \zeta_{s-1}(t)^{-1})^{-1}\cdots
(\zeta_1(t) v_2 \zeta_1(t)^{-1})^{-1} \mu(t)^{-1}
\end{align*}
By (1), $P_G(\mu\zeta_j)=B$ for all $j$,
and therefore each of the factors
$$\mu(t)(\zeta_{j}(t) v_{j+1} \zeta_{j}(t)^{-1})^{-1}\mu(t)^{-1}$$ admits a limit $1$.
So $\lambda(t) \mu(t) \lambda(t)^{-1}\mu(t)^{-1}$ admits a limit in $U$
by (3).

Next we consider the ``moreover'' part.
($\dagger$) holds for $g\in H$ as well.
So $P_H(\lambda) = u_s^{-1} P_H(\zeta_s) u_s$ is a parabolic subgroup of $H$.
\end{proof}

\begin{lem}
\label{lem:F-finiteness}
Let $F: M \to M$ be an automorphism of a connected reductive group.
Let $B\subset M$ be a Borel subgroup fixed by $F$, with unipotent radical $U$.
There exists a cocharacter $\mu$ of $M$, a positive integer $d$ and an element $u$ of $U$
such that $\mu = u F^d(\mu) u^{-1}$ and $B=P_M(\mu)$.
\end{lem}

\begin{proof}
By replacing $M$ by its derived subgroup, we can and do assume $M$ is semi-simple.
Let $\mu$ be a cocharacter of $M$
such that $B=P_M(\mu)$.

Let $i\ge 0$ be an integer.
There exists a maximal torus $T_i$ of $B$ such that $F^i(\mu)\subset T_i$.
Since all maximal tori of $B$ are conjugate by an element of $U$,
there exists an element $u_i$ of $U$ such that $T_0 = u_i T_i u_i^{-1}$.

So $u_i^{-1}F^i(\mu) u_i \subset T_0$,
and we can regard it as an element $x_i$ of the cocharacter lattice $X_*(M, T_0)$.
Since $\mu$ is a regular cocharacter, its centralizer $Z_M(\mu)$ is a maximal torus of $M$, and thus is just $T_0$.
Since automorphisms of $M$ send the centralizers to the centralizers,
$u_i^{-1} F^i u_i : M \to M$ fixes $T_0$.
Recall that $\Aut(M)\subset \operatorname{Inn}(M) \rtimes \Aut(\Dynkin(\Phi(M, T_0)))$, that is, after fixing a pinning,
an automorphism of $M$ comes from an automorphism of its Dynkin diagram.
Since $u_i^{-1} F^i u_i$ fixes $T_0$ and $B$,
it induces an isomorphism of the Dynkin diagram of $M$
and thus induces an isometry of the coroot lattice of $M$.
Since $M$ is semi-simple,
its coroot lattice and its cocharacter lattice span the same $\mathbb{R}$-vector space, and thus
$u_i^{-1} F^i u_i$ 
induces an isometry of $X_*(M, T_0)\otimes_{\mathbb{Z}}\mathbb{R}$.
In particular, the set $\{x_i\}$ is bounded and thus finite.
So $x_{i_0} = x_{i_0+d}$ for some $i_0 \ge 0$ and $d>0$.
We have $u_{i_0}^{-1}F^{i_0}(\mu) u_{i_0} = u_{i_0+d}^{-1}F^{i_0+d}(\mu) u_{i_0+d}$.
Thus $\mu = u_{i_0}u_{i_0+d}^{-1}F^d(\mu)u_{i_0+d}u_{i_0}^{-1}$.
\end{proof}

Recall a subgroup $\Gamma \subset G(\bar k)$ is said to be $G$-irreducible if
$\Gamma$ is not contained in any proper parabolic subgroup of $G(\bar k)$.

\begin{thm}
\label{thm:sw-irr}
Let $M$ be a connected reductive group over a field $k$.
Let $\bar k$ be the algebraic closure of $k$.
Let $F_M : M\to M$ be an automorphism of $M$ which can be realized
as conjugation by an element $g\in G(k)$ after an embedding
$M \hookrightarrow G$.

If $g$ and $Z_M(M)$ generate a $G$-irreducible subgroup,
then $M$ is a torus.
\end{thm}

\begin{proof}
One of the key ingredients is the results of Steinberg on endormorphisms of
linear algebraic groups.
By \cite[Theorem 7.2]{St68}, any automorphism of a linear algebraic group fixes a Borel subgroup.
Let $B_M\subset M$ be a Borel fixed by $F_M$.

There exists a cocharacter $\lambda:\Gm\to M$ such that $B_M = P_M(\lambda)$.
Let $U_M$ be the unipotent radical of $B_M$.
By the previous lemma, there exists $d>0$ and an element $u$ of $U_M$
such that $F_M^d(\lambda) = u \lambda u^{-1}$.
Consider the fake cocharacter $\mu:\Gm\to M$, defined by
$$
\mu := F_M^{d-1}(\lambda)F_M^{d-2}(\lambda) \cdots F_M(\lambda)\lambda.
$$
Note that
$F_M(\mu) = (u \lambda u^{-1}) \mu \lambda^{-1}$.

By Lemma \ref{lem:more-dynamic}, we have
\begin{itemize}
\item[(i)]  $P_G(\mu)$ is representable by a parabolic subgroup of $G$;
\item[(ii)] $P_M(\mu) = P_M(\lambda) = M \cap P_G(M)$.
\end{itemize}

\paragraph{\bf Claim} $g \in P_G(\mu)$.

\begin{proof}
We verify this using the definition of $P_G(\mu)$.
We have
\begin{align*}
\lim_{t\to 0}\mu(t) g  \mu(t)^{-1}
&=
\lim_{t\to 0}\mu(t) g  \mu(t)^{-1} g  ^{-1} g  \\
&=
\lim_{t\to 0}\mu(t)F_M(\mu)(t)^{-1} g  \\
&=
\lim_{t\to 0}\mu(t) \lambda(t) \mu(t)^{-1}u\lambda(t)^{-1}u^{-1}  g  \\
&=
\lim_{t\to 0}(\mu(t) \lambda(t) \mu(t)^{-1}\lambda(t)^{-1})(\lambda(t)u\lambda(t)^{-1})u^{-1}  g  \\
\end{align*}
The claim follows from Lemma \ref{lem:more-dynamic} (4).
\end{proof}

Note that since $\mu$ is valued in $M$, $Z_M(M) \subset Z_G(\mu)$.

Let $\Gamma$ be the subgroup of $G$ generated by
$Z_M(M)^\circ$ and $g$.
As a consequence of the claim, we have $\Gamma\subset P_G(\mu)$.
By Lemma \ref{lem:more-dynamic} (1), 
$P_M(\mu) = P_M(\lambda)$ is a Borel subgroup of $M$.
By the dynamic description of Borel subgroups, we have $P_G(\mu) \cap M = P_M(\mu)$.
So $P_G(\mu)$ is a proper parabolic subgroup of $G$
if $P_M(\mu)$ is a proper parabolic subgroup of $M$.
Since $\Gamma$ is assumed to be $G$-irreducible,
we must have $P_M(\mu) = M$.
Since $M=P_M(\mu)=B_M$ is chosen to be a Borel subgroup of $M$,
$M$ is forced to be a torus.
\end{proof}

\begin{cor}
\label{cor:sw-max-rank}
Let $M$ be a connected reductive group over a field $k$.
Let $\bar k$ be the algebraic closure of $k$.
Let $F_M : M\to M$ be an automorphism of $M$ which can be realized
as conjugation by an element $g\in G(k)$ after an embedding
$M \hookrightarrow G$.

Assume
\begin{itemize}
\item [(i)]
$g$ and $Z_M(M)^\circ$ generate a $G$-completely reducible subgroup;
\item [(ii)] $\rk M = \rk G$ and $\operatorname{char} k \ne 2,3$; and 
\item [(iii)] $M$ has a connected center.
\end{itemize}
Then $F_M$ fixes a maximal torus $T$ of $M_{\bar k}$.
\end{cor}

\begin{proof}
Let $\Gamma$ be the subgroup of $G$ generated by
$Z_M(M)$ and $g$.
If $\Gamma$ is $G$-irreducible, we are done because of Theorem \ref{thm:sw-irr}.
So we assume there exists a proper parabolic subgroup $P$ of $G_{\bar k}$ such that $\Gamma \subset P$.

By Borel-de Siebenthal theory (see \cite{Pep15} or \cite[Theorem 0.1]{Gil10}), when $\operatorname{k}\ne 2,3$, $\rk M = \rk G$ implies
$M=Z_G(Z_M(M))^\circ$.

We will prove a slightly stronger version of the corollary. We claim $F_M$ fixes a maximal torus of $M_{\bar k}$ assuming (i), (ii), and

\begin{quotation}
(iii') There exists a torus $Z$ of $M$ such that $M = Z_G(Z)^\circ$.
\end{quotation}

Since $\Gamma$ is $G$-completely reducible, there exists a Levi subgroup $L\subset P$
such that $\Gamma \subset L$.
Note that $(M \cap L)^\circ = Z_L(Z)^\circ$, which is a reductive subgroup (see \cite[Lemma 0.2(1)]{Gil10}) of $L$ fixed by $g$.
We claim $(M\cap L)^\circ$ is of maximal rank.
Let $S$ be any maximal torus of $L$ containing $Z$.
Since $S$ is commutative and connected, we have $S \subset Z_L(Z)^\circ = (M \cap L)^\circ$.
Thus $\rk (M \cap L)^\circ = \rk S =  \rk L = \rk G$.

We apply induction on the dimension of $G$.
Since $Z_L(Z)^\circ = (M\cap L)^\circ$, assumption (iii') is satisfied by $(M\cap L)^\circ$;
assumption (i) is also satisfied because~$L$ is a Levi subgroup of $G$.
Since $\dim L < \dim G$, by induction
there exists a maximal torus
$T$ of $(M \cap L)^\circ$ which is fixed by $F_M$.
Since $\rk (M \cap L)^\circ = \rk G$,~$T$ is also a maximal torus of~$G$.
\end{proof}

\subsection{}
We explain how our methods can possibly be used to establish a stronger form of
Steinberg-Winter, at least for groups having connected center.
Dynamic methods are very well behaved for disconnected linear algebraic groups.
We similarly define $G$-complete reducibility for general linear algebraic groups
by replacing parabolics by pseudo-parabolics.
Let $F: M\to M$ be an automorphism
which can be realized as conjugation by an element $g$ of $G$ after an embedding $M \hookrightarrow G$.
Let $H$ be the scheme-theoretic closure of the (abstract) group generated by $M$ and $g$.
Note that $H$ is a disconnected reductive group, and $\rk H = \rk M$.
Let $\Gamma$ be the subgroup of $H$ generated by $Z_M(M)^\circ$ and~$g$.
We expect that the $H$-complete reducibility of $\Gamma$ implies the existence of a fixed maximal torus.

\section{The structure of $G$-completely reducible mod $\varpi$ Galois representations}

In this section, we give a complete description of all 
$G$-completely reducible mod $\varpi$ Galois representations
valued in split reductive groups.

The first step is to show $G$-complete reducibility implies tame ramification,
reducing the classification of mod $\varpi$ Galois representations
to the question of classification of (certain) solvable subgroups of derived length $2$ of reductive groups.

\begin{lem}
\label{lem:trivial-wild-inertia}
Let $P_K$ be the wild inertia of $G_K$.
If $\bar\rho:G_K\to G(\bFp)$ is $G$-completely reducible, $\bar\rho(P_K) = \{\Id\}$.
\end{lem}

\begin{proof}
Let $P_K\subset G_K$ be the wild inertia.
The image $\bar\rho(P_K) \subset G(\bFp)$ is a $p$-group, and thus consists of unipotent elements.
By \cite[Corollaire 3.9]{BT71}, there exists a parabolic subgroup $P$ of $G_{\bFp}$ with unipotent radical $R_u(P)$
such that
\begin{itemize}
\item $\bar\rho(P_K)\subset R_u(P)(\bFp)$, and
\item $N(\bar\rho(P_K))\subset P(\bFp)$;
\end{itemize}
here $N(\bar\rho(P_K))$ is the normalizer of $\bar\rho(P_K)$.
Since $P_K$ is a normal subgroup of $G_K$, $\bar\rho(G_K)\subset N(\bar\rho(P_K))\subset P(\bFp)$.
Since $\bar\rho$ is $G$-completely reducible, $\bar\rho(G_K)$ is contained in a Levi subgroup $L$ of $P$.
So $\bar\rho(P_K) \subset L(\bar\bF_p) \cap R_u(P)(\bFp) = \{\Id\}$.
\end{proof}

\begin{dfn}
\label{dfn:strongly-semisimple}
We say $\bar\rho: G_K \to G(\bFp)$ is \textit{quasi-semisimple}
if there exists a maximal torus $T$ of $G(\bFp)$ such that
$\bar\rho(I_K) \subset T(\bFp)$ and
$\bar\rho(G_K) \subset N_G(T(\bFp))$.
\end{dfn}

\begin{thm}
\label{thm:strongly-semisimple}
If $\bar\rho: G_K \to G(\bFp)$
is $G$-completely reducible, 
then $\bar\rho$ is quasi-semisimple.

Moreover, if $\bar\rho$ is $G$-irreducible,
there exists a \textit{unique} maximal torus $T$ of $G(\bFp)$ containing $\bar\rho(I_K)$.
Consequently, if $\bar\rho(G_K)\subset G(\bF)$,
$T$ has a model defined over the ring of Witt vectors $W(\bF)$.
\end{thm}

\begin{proof}
By induction on the dimension of $G$,
we can reduce the general case to the case where $\bar\rho$ is $G$-irreducible.
Recall that $\bar\rho$ is $G$-irreducible if it does not factor through any proper parabolic of $G$.
If $\bar\rho$ does factor through a proper parabolic of $G$, the $G$-complete reducibility forces
$\bar\rho$ to factor through a proper Levi subgroup of $G$, which is a reductive group
of strictly smaller dimension.

So we assume $\bar\rho$ is $G$-irreducible in the rest of the proof.
By Lemma \ref{lem:trivial-wild-inertia}, $\bar \rho(I_K)$ is a finite cyclic group
generated by elements of order prime to $p$.
Write $M$ for $Z^{\circ}_{G(\bFp)}(\bar \rho(I_K))$,
the neutral component of the centralizer of $\bar \rho(I_K)$ in $G$.
Since $\bar\rho(I_K)$ consists of semi-simple elements of $G(\bF)$,
$M$ is a reductive subgroup of~$G$.
Let $\Phi_K \in G_K$ be a topological generator of $G_K/I_K$.
Since $I_K$ is a normal subgroup of~$G_K$,
the conjugation by $\bar\rho(\Phi_K)$ action induces an automorphism of~$M$,
which we denote by $F_M: M \to M$.

Next we show $\bar\rho(I_K)\subset Z_M(M)$.
Since $G$ is connected, a semisimple element of~$G$ is contained in a maximal torus.
Since $\bar\rho(I_K)$ is a cyclic group consisting of semi-simple elements,
there exists a maximal torus $T$ containing $\bar\rho(I_K)$.
Since a torus is connected, we have $T\subset M$, and thus $\bar\rho(I_K)\subset M(\bFp)$.
It is immediate from the definition of $M$ that $\bar\rho(I_K) \subset Z_M(M)$.

By Theorem \ref{thm:sw-irr}, $M$ is a torus.
Let $T$ be any maximal torus of $G$ containing $\bar\rho(I_K)$.
Since $T$ is commutative and connected, we have $T \subset Z_G(\bar\rho(I_K))^\circ = M$.
So $M$ is the unique maximal torus containing $\bar\rho(I_K)$.
Now consider the ``moreover'' part.
For $\sigma \in \Gal(\bFp/\bF)$,
$\sigma(M)$ is also a maximal torus containing $\bar\rho(I_K)$.
So $\sigma(M) = M$, and thus by Galois descent $M$ is defined over $\bF$.
By \cite[B.3.5]{Crd14}, $T$ has a model over $W(\bF)$.
\end{proof}

\subsection{Example}
We illustrate the technical proof using a very concrete example.
Let $G=\GL_4$.
Let $\bar\rho:G_K \to \GL_4(\bFp)$ be a semi-simple Galois representation.
We decompose $V=V_{\chi_1}\oplus V_{\chi_2}$ into $I_K$-isotropic subspaces.
Here $\chi_1$, $\chi_2:I_K \to \bFp^\times$ are distinct charatcers 
such that for $v\in V_{\chi_i}$ and $\sigma\in I_K$,
$\bar\rho(\sigma)v = \chi_i(\sigma)v$, $i=1,2$.
$$
\bar\rho|_{I_K} = 
\begin{bmatrix}
\chi_1 & & & \\
& \chi_1 & & \\
& & \chi_2 & \\
& & & \chi_2
\end{bmatrix}
$$
There are two possibilities: either both $V_i$ are $\bar\rho(\Phi_K)$-stable, or
$\bar\rho(\Phi_K)$ sends $V_i$ to $V_{3-i}$, $i=1,2$.
The first case is simple: $V=V_{\chi_1}\oplus V_{\chi_2}$ as a $G_K$-module. 
Now we consider the latter case.
By Steinberg's theorem \cite[Theorem 7.2]{St68}, we can assume
$\bar\rho(\Phi_K)$ fixes a Borel
$$
P_M = \begin{bmatrix}
* & * & & \\
& * & & \\
& & * & * \\
& & & *
\end{bmatrix}
$$
of $M=\GL_2 \times \GL_2 \subset \GL_4$ and thus we must have
$$
\bar\rho(\Phi_K) = 
\begin{bmatrix}
 & & a & b \\
 & & & c \\
 d & f & & \\
 & e & &
\end{bmatrix}
$$
for some $a,b,c,d,e,f\in \bFp$.
The Borel subgroup $P_M$ is of shape $P_M(\lambda)$
for 
$$
\lambda(t) = \begin{bmatrix}
t^{\alpha} & * & & \\
& t^{\beta} & & \\
& & t^{\gamma} & * \\
& & & t^{\delta}
\end{bmatrix}
$$
for $\alpha > \beta$ and $\gamma > \delta$.
We have
$$
\bar\rho(\Phi_K)\lambda(t)\bar\rho(\Phi_K)^{-1} = \begin{bmatrix}
t^{\gamma} & * & & \\
& t^{\delta} & & \\
& & t^{\alpha} & * \\
& & & t^{\beta}
\end{bmatrix}
$$
and thus
$$
\bar\rho(\Phi_K)\lambda(t)\bar\rho(\Phi_K)^{-1}\lambda(t) = \begin{bmatrix}
t^{\alpha+\gamma} & * & & \\
& t^{\beta+\delta} & & \\
& & t^{\alpha+\gamma} & * \\
& & & t^{\beta+\delta}
\end{bmatrix}
$$
Since $\alpha+\gamma>\beta+\delta$, we have
$$
P_{\GL_4}(\bar\rho(\Phi_K)\lambda\bar\rho(\Phi_K)^{-1}\lambda) = 
\begin{bmatrix}
* & * & * & *\\
& * & & * \\
* & * & * & * \\
& * && *
\end{bmatrix}
$$
and finally we observe $\bar\rho(\Phi_K) \in P_{\GL_4}(\bar\rho(\Phi_K)\lambda\bar\rho(\Phi_K)^{-1}\lambda)$.
In general, by Lemma \ref{lem:F-finiteness}, there exists an integer $d$ such that
$\prod_{i=d-1}^0 \bar\rho(\Phi_K)^i\lambda(t)\bar\rho(\Phi_K)^{-i}$ gives the desired parabolic.

\section{Crystalline lifts of irreducible mod $\varpi$ Galois representations}

\label{subsec:crys-ch}

Write $\kappa$ for the residue field of $K$.
Fix a coefficient field $E$ with ring of integers~$\cO$ and uniformizer $\varpi$.
Write $\bF$ for the residue field $\cO/\varpi$.
Assume $\kappa\subset \bF$.
Let $\Phi_K \in G_K$ be a (lift of a) topological generator of $G_K/I_K$.
Fix an algebraic closure $\bar K$ of $K$.

In this section, we assume $G$ is a split group since we are primarily interested in Galois representations valued in $L$-groups.
The $L$-group of a connected reductive group is split, albeit possibly disconnected.

\subsection{}
Let $T$ be a maximal torus of $G$. More precisely, $T$ is a smooth group scheme over $\spec \cO$ such that $T_{\bar k}\subset G_{\bar k}$ is a maximal torus for all geometric points $\cO \to \bar k$.
Write
$W(G, T)$ for the Weyl group scheme $N_G(T)/T$.
Write $\mathcal{O}_{\bar K}$ for the ring of integers in $\bar K$.
Note that we have a commutative diagram
$$
\xymatrix{
    W(G, T)(\mathcal{O})
    \ar@{^{(}->}[r] \ar[d] &
    W(G, T)(\mathcal{O}_{\bar K})\ar@{=}[d] \ar@{=}[r] & W(G, T)({\bar K}) \\
    W(G, T)(\bF)
    \ar@{^{(}->}[r] &
    W(G, T)(\bFp) \\
}
$$
and as a consequence, the map
$W(G, T)(\mathcal{O}) \to W(G, T)(\bF)$ is injective.
On the other hand, since $N_G(T)$ is a smooth group scheme over $\spec \cO$, 
$W(G, T)(\mathcal{O}) \to W(G, T)(\bF)$ is also surjective.
We will identify $W(G, T)(\mathcal{O})$ with $W(G, T)(\bF)$
and write it as $W(G, T)$.
It will be clear from the context if the notation $W(G, T)$ denotes a set or a group scheme.

Write $M_{T,\cris}$ for the set of representations $I_K \to T(\cO)$,
which can be extended to a crystalline representation $G_{K'}\to T(\cO)$
for some finite unramified extension $K'/K$ inside $\bar K$.
Since the union of two finite unramified extensions inside $\bar K$ is still a finite unramified extension,
$M_{T,\cris}$ is an abelian group.

The abelian group $M_{T,\cris}$ has a $\bZ[W(G,T)]$-module structure, defined by
$w v := (\sigma \mapsto w v(\sigma) w^{-1})$,
for $w\in W(G,T)$ and $v \in M_{T,\cris}$.

The abelian group $M_{T,\cris}$ also has a $\bZ[G_K/I_K]$-module structure, defined by
$\alpha v := (\sigma \mapsto v(\alpha^{-1} \sigma \alpha))$
for $\alpha\in G_K$ and $v \in M_{T,\cris}$.

The following lemma is clear.

\paragraph{\bf Lemma and Definition}
The $\bZ[W(G,T)]$-structure and the $\bZ[G_K/I_K]$-structure
on $M_{T,\cris}$ commute with each other.
Therefore $M_{T,\cris}$ is a $\bZ[W(G,T)] \otimes \bZ[G_K/I_K]$-module.

Similarly,
write $M_{T,\modp}$ for the abelian group of mod $\varpi$ representations $I_K \to T(\bF)$.
The abelian group $M_{T, \modp}$ has a $\bZ[W(G,T)] \otimes \bZ[G_K/I_K]$-module structure.

\begin{lem}
\label{lem:extendable}
Write $\zeta: N_G(T) \to W(G,T)$ for the quotient map.

(1)
Let $w$ be an element of $N_G(T)(\cO)$ of finite order.
An element $v\in M_{T,\cris}$ extends to
a continuous representation $\rho: G_K \to N_G(T)(\cO)$
by setting $\rho(\Phi_K)=w^{-1}$ and $\rho|_{I_K} = v$
if and only if
$$v \in \ker(M_{T,\cris} \xrightarrow{\zeta(w) \otimes 1 - 1 \otimes \Phi_K} M_{T,\cris}).$$

(2)
Let $\bar w$ be an element of $N_G(T)(\bF)$.
An element $v\in M_{T,\modp}$ extends to
a representation $\bar\rho: G_K \to N_G(T)(\bF)$
by setting $\bar\rho(\Phi_K)=\bar w^{-1}$ and $\bar\rho|_{I_K} = v$
if and only if
$$v \in \ker(M_{T,\modp} \xrightarrow{\zeta(\bar w) \otimes 1 - 1 \otimes \Phi_K} M_{T,\modp}).$$
\end{lem}

\begin{proof}
(1)
Since $w$ is of finite order,
it suffices to show
$v\in M_{T,\cris}$ extends to
a representation $\rho: W_K \to N_G(T)(\cO)$ of the Weil group $W_K \cong I_K \rtimes \bZ$
by setting $\rho(\Phi_K)=w^{-1}$ and $\rho|_{I_K} = v$
if and only if
$v \in \ker(M_{T,\cris} \xrightarrow{\zeta(w) \otimes 1 - 1 \otimes \Phi_K} M_{T,\cris})$.
If $v$ is extendable to $\rho$, then for all $\sigma\in I_K$
$$\rho(\Phi_K^{-1} \sigma \Phi_K) = w \rho(\sigma) w^{-1};$$
the left hand side restricted to $I_K$ is $(1\otimes \Phi_K) v$,
and the right hand side restricted to $I_K$ is $(\zeta(w) \otimes 1) v$.
So $(1\otimes \Phi_K) v = (\zeta(w) \otimes 1) v$.
Conversely, if 
$(1\otimes \Phi_K) v = (\zeta(w) \otimes 1) v$, then
$v(\Phi_K^{-1} \sigma \Phi_K) = w v(\sigma) w^{-1}$
for all $\sigma\in I_K$.
Define $\rho(\sigma \Phi^n):=v(\sigma)w^{-n}$ for all $\sigma\in I_K$ and $n\in \bZ$.
It is clear $\rho$ is well-defined on $W_F$, and extends to $G_F$ uniquely by continuity.

(2) is similar to (1).
\end{proof}

\begin{dfn}
For an element of the Weyl group $\mathpzc{w} \in W(G,T) = W(G,T)(\bF) = W(G,T)(\cO)$,
define 
$$
M_{T, \mathpzc{w}, \cris} := \ker(M_{T,\cris} \xrightarrow{\mathpzc{w} \otimes 1 - 1 \otimes \Phi_K} M_{T,\cris})
\text{, and}
$$
$$
M_{T, \mathpzc{w}, \modp} := \ker(M_{T,\modp} \xrightarrow{\mathpzc{w}\otimes 1 - 1 \otimes \Phi_K} M_{T,\modp}).
$$
\end{dfn}

The following simple lemma is essentially how we construct crystalline lifts.

\begin{lem}
\label{lem:simple-trick}
Let $\bZ[X]$ be the polynomial ring.
Let $a(X),b(X)\in \bZ[X]$ be two polynomials.
Let $n$ and $N$ be integers.
Assume $a(n)b(n) = 0$.

Let $\widetilde{M}$ be a $\bZ[X]/(a(X)b(X) - N)$-module.
Write $M$ for $\widetilde{M}\otimes_{\bZ}\bZ/N$.

If $\widetilde{M}$ has a torsion-free, finitely generated underlying abelian group,
the sequence
$$
0 \to a(X)M \to M \xrightarrow{\cdot b(X)} b(X) M \to 0
$$
is short exact.
\end{lem}

\begin{proof}
Pick $\bar v \in \ker (M \to b(X) M)$.
Let $v \in \widetilde{M}$ be a lifting of $\bar v$.
We have $b(X)v \mapsto 0$ in $M$.
Since $M = \widetilde{M} \otimes \bZ/N$,
$b(X)v = N u$ for some $u\in \widetilde{M}$.
Multiply both sides by $a(X)$, we get
$a(X)b(X)v = Nv = N a(X) u$.
Since $\widetilde{M}$ is $\bZ$-torsion-free, we have $v = a(X)u$, as desired.
\end{proof}

\begin{prop}
\label{prop:lift-inertia}
If $\mathpzc{w}^{[\bF: \kappa]}=1$ and $E$ contains $K$,
the map
$$
M_{T, \mathpzc{w}, \cris} \to M_{T, \mathpzc{w}, \modp}
$$
is surjective.
\end{prop}

\begin{proof}
Write $f:= [\bF:\kappa]$. Let $K_f$ be the unramified extension of $K$ of degree~$f$.

We single out a $\bZ[W(G,T)] \otimes \bZ[G_K/I_K]$-submodule $M_{T, \cris}^{0} \subset M_{T, \cris}$
which consists
of elements that can be extended to a representation $G_{K_f}\to T(\cO)$.
Note that
$
M_{T, \cris}^{0} \to M_{T, \modp}
$
is surjective because the fundamental character of niveau $f$ 
admits a crystalline lift, namely, the Lubin-Tate character of the field~$K_f$.
Put $M_{T, \mathpzc{w}, \cris}^0 := M_{T, \cris}^{0} \cap M_{T, \mathpzc{w}, \cris}$.

Note that on both $M_{T,\cris}^0$ and
$M_{T,\modp}$, we have $(\mathpzc{w} \otimes 1)^f = (1 \otimes \Phi_K)^f = \Id$,
where $\Phi_K$ is the fixed topological generator of $G_K/I_K$.

Put
$$
\Xi := \sum_{i = 0}^{f-1} \mathpzc{w}^i \otimes \Phi_K^{f-1-i}.
$$
Commutativity of $\mathpzc{w}\otimes 1$ and $1\otimes \Phi_K$
implies $(\mathpzc{w} \otimes 1 - 1 \otimes \Phi_K)\Xi = (\mathpzc{w} \otimes 1)^f - (1 \otimes \Phi_K)^f$.
In particular, the inclusion 
$\Xi M_{T,\cris}^0 \to M_{T,\cris}^0$
factors through $M_{T, \mathpzc{w}, \cris}^0$ (which is the arrow at the top of the diagram below).

Consider the commutative diagram
$$
\xymatrix{
\Xi M_{T,\cris}^0 \ar[r] \ar[d] &
M_{T, \mathpzc{w}, \cris}^0 \ar[d]\ar[r] &
M_{T,\cris}^0\ar[d] \ar[rr]^{\mathpzc{w} \otimes 1 - 1 \otimes \Phi_K} && M_{T,\cris}^0 \ar[d]\\
\Xi M_{T,\modp} \ar[r]\ar[d] &
M_{T, \mathpzc{w}, \modp} \ar[r] &
M_{T,\modp} \ar[rr]^{\mathpzc{w} \otimes 1 - 1 \otimes \Phi_K} && M_{T,\modp} \\
0 & & &
}
$$
It is clear that $\Xi M_{T,\cris}^0 \to \Xi M_{T, \modp}$ is surjective.
So it suffices to show
$$
\Xi M_{T,\modp} \hookrightarrow
M_{T, \mathpzc{w}, \modp}
$$
is surjective.

Let $\bar \chi: I_K \to \bF^\times$
be a fundamental character of niveau $f$.
Note that $\bar\chi$ generates the abelian group $M_{\Gm, \bF}$.
Indeed, there is an abelian group isomorphism
$\iota_{\bar\chi}:\bZ/(q^f-1) \xrightarrow{\cong} M_{\Gm, \bF}$,
sending $1$ to $\bar\chi$.
We have $M_{T, \bF} \cong M_{\Gm, \bF} \otimes_{\bZ} \Hom_{\text{GrpSch}}(\Gm, T)$.
Note that the Weyl group element $\mathpzc{w}$ acts on $\Hom_{\text{GrpSch}}(\Gm, T)$ via conjugation $v \mapsto \mathpzc{w} v \mathpzc{w}^{-1}$.

We specialize Lemma \ref{lem:simple-trick} as follows:
\begin{itemize}
\item Set $\widetilde{M} = \Hom_{\text{GrpSch}}(\Gm, T)$, and regard it as a $\bZ[X]$-module where $X$ acts by $\mathpzc{w}$;
\item Set $M = M_{T,\modp}$, and regard $M$ as a $\bZ[X]$-module via $X \mapsto \mathpzc{w} \otimes 1$;
\item Set $N = q^f-1$;
\item Set $n = q$;
\item Set $a(X) = \sum_{i=0}^{f-1-i}X^iq^{f-1-i}$;
\item Set $b(X) = q-X$;
\end{itemize}
We can identify $M$ with $\widetilde{M} \otimes_{\bZ} \bZ/(q^f-1)$ via
the map $\iota_{\bar\chi}:\bZ/(q^f-1) \xrightarrow{\cong} M_{\Gm, \bF}$.

Here are a few things to check:
\begin{itemize}
\item [(i)] $\widetilde{M}$ is finitely generated and torsion-free over $\bZ$.
\item [(ii)] $(a(X)b(X)- q^f+1)$ kills $\widetilde{M}$;
\item [(iii)] $\widetilde{M} \otimes_{\bZ}\bZ/(q^f-1) \cong M$ as abelian groups;
\item [(iv)] $a(q)b(q) = 0$.
\end{itemize}
Items (i), (iii) and (iv) are clear.
For item (ii), notice that
$a(X)b(X)=q^f-X^f$. Since we assumed $\mathpzc{w}^f=1$, $a(X)b(X)=q^f-1$.
\end{proof}

The goal of the rest of this section is to prove the following theorem:

\begin{thm}
\label{thm:irreducible-lift}
Let $\kappa$ be the residue field of $K$.
Let $\bF/\kappa$ be a finite extension.
Let $K^{\ur}$ be the maximal unramified extension of $K$ with ring of integers $\cO_{K_{\ur}}$.

Let $\bar\rho : G_K \to G(\bF)$ be a quasi-semisimple (see Definition \ref{dfn:strongly-semisimple}) representation.

(1)
There exists a crystalline representation $\rho: G_K \to G(\cO_{K^{\ur}})$
lifting $\bar\rho$.

(2)
Assume $G$ admits a simply-connected derived subgroup and $\bar\rho$ is $G$-irreducible.
Let $\bF_{\bar\rho}$ be the splitting field of
$\bar\rho|_{I_K}$, that is, the smallest
field extension $\bF_{\bar\rho}$ of $\mathbb{F}$ such that $\bar\rho|_{I_K}: I_K \to G(\bF)$ factors through the $\bF_{\bar\rho}$-points of a split torus of $G$.
Then $\rho$ can be chosen to
have image in $G(\cO_{K_{\bar\rho}})$
where $K_{\bar\rho}$ is the unramified extension of $K$ with residue field $\bF_{\bar\rho}$.
\end{thm}

\subsection{}

The strategy is as follows: 
the first step is to choose a lift of $\bar\rho|_{I_K}$
which admits an extension to the whole Galois group $G_K$.
This is already done in Proposition~\ref{prop:lift-inertia}.
The second step is to choose a lift of all Frobenius elements.
The continuity of the lift is free because we'll
only use
finite order lifts (modulo the image of~$I_K$) of Frobenius elements.

\begin{lem}
\label{lem:extended-ng}
Assume the special fiber $T_{\bF}$ of $T$
is a split torus.
There exists a finite subgroup $\widetilde{N} \subset N_G(T)(W(\bF))$
such that $\widetilde{N} \to N_G(T)(\bF)$ is surjective.
\end{lem}

\begin{proof}
By \cite[B.3.5]{Crd14}, $T$ splits if and only if $T_{\bF}$ splits.
The key ingredient is Tits' theory of
extended Weyl groups.

By \cite{Ti66}, there exists a subgroup $\widetilde{W}\subset N_G(T)(W(\bF))$ which is isomorphic to
the extension of the Weyl group $W(G,T)$ by $(\bZ/2)^{\otimes l}$ for some $l\ge 0$,
and generates the whole Weyl group.
Write $[-]: T(\bF) \to T(W(\bF))$ for the Teich\"muller lift.

\paragraph{\bf Fact}
The Teichm\"uller lift is the unique $p$-adic continuous multiplicative section of $T(W(\bF))\to T(\bF)$.

\begin{proof}
We include a proof here because it is short.
It is well-known for $T=\Gm$. 

In general, choose a faithful representation $i:T \to \GL_N \subset \Mat_{N\times N}$.
Let $s,t: T(\bF) \to T(\cO)$ be two multiplicative sections.
We have $i(s(x)) - i(t(x)) \equiv 1 \mod p^f$ for all $x\in  T(\bF)$;
$(i(s(x)) - i(t(x)))^{p^{nf}} \equiv 1 \mod p^{(n+1)f}$; and
$i(s(x)) - i(t(x)) = i(s(x^{p^{nf}})) - i(t(x^{p^{nf}})) \equiv (i(s(x)) - i(t(x)))^{p^{nf}} \equiv 1 \mod p^{nf}$ for all~$n$.
\end{proof}

For each $w \in \widetilde{W}$ and $x\in T(\bF)$,
$x\mapsto w^{-1}[wxw^{-1}]w$ is a continuous section of $T(W(\bF))\to T(\bF)$
and must be equal to the Techm\"uller lift.
Let $\widetilde{N}$ be the composite $\widetilde{W}\cdot [T(\bF)]$.
Since for all $w,w'\in \widetilde{W}$ and $x,x'\in T(\bF)$, we have
$w[x]w'[x'] = ww'[w'^{-1}xw'x']$, $\widetilde{N}$ is a finite order subgroup of
$N_G(T)(W(\bF))$, as desired.
\end{proof}

The existence of $\widetilde{N}$ has the following immediate consequence:

\begin{cor}
\label{cor:lift-witt}
Let $\bar\rho: G_K \to G(\bF)$ be a $G$-completely reducible representation.
There exists a lift $\rho: G_K \to G(W(\bFp))$ of $\bar\rho$.

Indeed, for any lift $v$ of $\bar\rho|_{I_K}$
to $G(\cO_{K^{\ur}})$ that can be extended to the whole Galois group $G_K$,
there exists a lift $\bar\rho$ to $G(\cO_{K^{\ur}})$
whose inertia is $v$.
\end{cor}

\begin{proof}
We first prove the first paragraph.
We are allowed to enlarge the coefficient field $\bF$ to make $T$ split.
Set $\rho|_{I_K}$ to be the Teichm\"uller lift of $\bar\rho|_{I_K}$.
Let $\Phi_K\in G_K$ be a lift of the topological generator of $G_K/I_K$.
Choose an element $n \in \widetilde{N}$ which lifts $\bar\rho(\Phi_K)$.
Set $\rho(\Phi_K) = n$.
Write $n=wt$ where $w$ is an element of Tits' extended Weyl group $\widetilde{W}$
and $t$ lies in the Teichm\"uller lift of $T(\bF)$.
Let $\sigma$ be an element of $I_K$.
Write $x$ for $\bar\rho(\sigma)$.
We have $\rho(\Phi_K\sigma \Phi_K^{-1}) = [\bar\rho (\Phi_K\sigma \Phi_K^{-1})]
=[w x w^{-1}]$.
By the proof of the previous lemma,
$[w x w^{-1}] = w [x] w^{-1} = w \rho(\sigma) w^{-1} = n \rho(\sigma) n^{-1}$,
and thus $\rho$ extends uniquely to a continuous homomorphism $G_K \to G(W(\bF))$.

Now we prove the ``indeed'' part.
It is an immediate consequence of Lemma~\ref{lem:extendable} and Lemma~\ref{lem:extended-ng}.
\end{proof}

\begin{lem}
\label{lem:exp-weyl}
Let $\bar\rho:G_K \to G(\bF)$ be a $G$-irreducible Galois representaion.
By Theorem \ref{thm:strongly-semisimple}, there exists a unique maximal torus $T$ of $G$ such that $\bar\rho(G_K)\subset N_G(T)(\bF)$.

Let $\kappa$ be the residue field of $K$.
Let $\bF_0\subset \bF$ be the smallest subfield of $\bF$ containing $\kappa$ 
such that
$\bar\rho(I_K) \subset T(\bF_0)$.
(Recall that $G$ is a Chevalley group and has a $\bZ$-model.)
Let $\Phi_K\in G_K$ be a lift of a topological generator of $G_K/I_K$.
The map $G_K \to N_G(T)(\bF) \to W(G,T)(\bF)$ maps
$\Phi_K$ to an element $w$ of the Weyl group $W(G,T)(\bF)$.

If $G$ admits a simply-connected derived subgroup, then $w^{[\bF_0:\kappa]}=1$ in $W(G,T)(\bF)$.
\end{lem}

\begin{proof}
Write $f_0 := [\bF_0:\kappa]$.
Let $s\in \bar\rho(I_K)$ be a generator.
By the proof of Theorem \ref{thm:strongly-semisimple},
$T = Z_G(s)^{\circ}$ is the connected centralizer of $s$.
Since $\bar\rho(I_K) \subset T(\bF_0)$,
we have $\bar\rho(\Phi_K)^{f_0}s\bar\rho(\Phi_K)^{-f_0} = s$.
So $\bar\rho(\Phi_K)^{f_0} \in Z_G(s) \cap N_G(T)$.
Since~$G$ has a simply-connected derived subgroup,
$Z_G(s) = Z_G(s)^\circ$.
So $\bar\rho(\Phi_K)^{f_0} \in T$, that is,
$w^{[\bF_0:\kappa]}=1$ in $W(G,T)(\bF)$.
\end{proof}

\begin{proof} [Proof of Theorem \ref{thm:irreducible-lift}]
(1)
We choose a sufficiently large coefficient field $E$ (which is unramified over $K$) such that the cardinality of the Weyl group $W(G,T)$ divides $[\bF:\kappa]$ and $T_{\bF}$ splits. 
The assumption of Proposition \ref{prop:lift-inertia} is satisfied.
So there exists a crystalline lift 
$v: I_K \to T(\cO)$ such that $v = \bar\rho|_{I_K}$ mod $\varpi$.
By Lemma \ref{lem:extendable}, $v$ can be extended to $G_K$.

(2)
For ease of notation, replace $\mathbb{F}$
by $\bF_{\bar\rho}$.
Write $\cO$ for $\cO_{K_{\bar\rho}}$.
We choose the field $\bF_0$ as in Lemma \ref{lem:exp-weyl}.
Note that the maximal torus in Lemma \ref{lem:exp-weyl} is split:
let $S$ be a maximal split torus over $\bF$
such that $\bar\rho(I_K)\subset S(\bF)$;
since $T=Z_G(\bar\rho(I_K))^\circ$, we have
$T \supset S$; now since $G$ is a split group, we must have $S=T$.
Let $K_{f_0}$ be the unramified extension of $K$ of degree $[K_{f_0}:K]=[\bF_0:\kappa]$.
Let $\cO_0$ be the ring of integers of $K_{f_0}$.
We have $\cO_0\subset \cO$.
By the previous Lemma, 
Proposition~\ref{prop:lift-inertia} is applicable, and thus
there exists a lift $v:I_K \to T(\cO_0)$ such that $v = \bar\rho|_{I_K}$ mod $\varpi$
and $v$ admits an extension to a representation $G_K \to N_G(T)(\cO_0)$.
By Corollary~\ref{cor:lift-witt},
$v$ admits an extension to $G_K$ which lifts $\bar\rho$.

Fix $\Phi_K \in G_K$, a lift of a topological generator of $G_K/I_K$.
By Lemma~\ref{lem:extended-ng}, we choose a finite order lift
$X\in \widetilde{N}\subset N_G(T)(\cO)$ of $\bar\rho(\Phi_K)$.
Since the Weyl group scheme is a constant group scheme,
any two lifts of $\bar\rho(\Phi_K)$ have the same conjugation action on the
maximal torus $N_G(T)(\cO)$, and therefore
we can extend $v$ to a representation $G_K \to N_G(T)(\cO)$ by setting $\Phi_K \mapsto X$.
\end{proof}

\section{Hodge-Tate theory for Galois representations valued in reductive groups}

\label{sec:hodge-tate-cocharacter}

The Hodge-Tate theory for $\GL_N$ is reviewed in Appendix \ref{Hodge-Tate}.
In this section, we discuss Hodge-Tate theory for general reductive groups, and
show $G$-irreducible mod $\varpi$ Galois representations admit Hodge-Tate regular
crystalline lifts.

\subsection{First properties of Hodge-Tate cocharacters}

\begin{dfn}
\label{def:ht}
Fix an algebraic closure $\overline{\Qp}$ of $\Qp$.
Let $K, E \subset \overline{\Qp}$ be finite extensions of $\Qp$.
The field $E$ will serve as the coefficient field.
To define colabeled Hodge-Tate gradings,
we assume $K$ is a subfield of $E$ and therefore $G_E$ as a subgroup of~$G_K$.

Let $\bC := \bC_K$ be the completed algebraic closure of $K$.
Let $\sigma: E\hookrightarrow \bC$ be an embedding.
Let $(\rho, V)$ be a Hodge-Tate representation of $G_K$.
Then one can define \textit{the $\sigma$-colabeled Hodge-Tate grading} on $\bC\otimes_{\sigma,E}V$ 
by setting the $i$-th graded piece to be
$$
\Img((\bC(i)\otimes_{\sigma,E}V)^{G_E}\otimes_E \bC(-i) \to \bC\otimes_{\sigma,E}V)
$$
which is compatible with tensor product and duality.

Let $G$ be a reductive group over $E$.
A $G$-valued representation is \text{Hodge-Tate} if
for all representations $G\to\GL(V)$,
$V$ is a Hodge-Tate $G_K$-module.
Let $\rho: G_K\to G(E)$ be a Hodge-Tate $G$-valued representation.
Consider $G(\sigma)\circ\rho: G_K \to G(\bC)$.
By Tannakian theory, 
there is a cocharacter $\hT(\rho)^{\sigma}: \Gm\to G_{\bC}$,
such that for any faithful representation
$i:G\to \GL_N$, the composition
$i(\hT(\rho)^{\sigma})$ recovers the
Hodge-Tate grading on
$i(G(\sigma)\circ \rho):G_K\to\GL_N(\bC)$.

Set $\hT(\rho):= (\hT(\rho)^{\sigma})_{\sigma:E\hookrightarrow \bC}\in\prod_{E\hookrightarrow\bC}X_*(G_{\bC})$.
We call $\hT({\rho})$ the \textit{co-labeled Hodge-Tate cocharacter} of $\rho$. 
\end{dfn}

The formation of co-labeled Hodge-Tate cocharacters is clearly functorial in $G$.

\begin{lem}
\label{lem:ht-functorial}
Let $f:G\to H$ be a morphism of reductive groups over $E$.
If $\rho:G_K \to G(E)$ is a Hodge-Tate representation,
we have $\hT(f\circ \rho) = f(\hT(\rho))$.
\end{lem}

\begin{proof}
It follows immediately from Tannakian theory.
\end{proof}

\subsubsection{Regular cocharacter}
\label{HT-regular}
Let $H$ be a reductive group with maximal torus $S$.
A cocharacter $x\in X_{*}(H,S)$ is said to be \textit{regular}
if it is not killed by any root of $H$ (with respect to $S$).

We say $\rho$ is \textit{Hodge-Tate regular}
if for all $\sigma:E\hookrightarrow \bC$,
the cocharacter $\hT(\rho)^{\sigma}$ of $G_{\bC}$ is regular.

When $G=\GL_N$, we can also define \textit{labeled Hodge-Tate weights} (see Appendix~\ref{Hodge-Tate}).
It turns out labeled Hodge-Tate regularity is equivalent to colabeled Hodge-Tate regularity.
So our definition coincides with the usual notion of Hodge-Tate regularity
in the literature.

\begin{lem}
Assume $G = \GL_N$.
Assume $E$ admits an embedding of the Galois closure of~$K$.
Then $\rho$ is Hodge-Tate regular if and only if 
the \textit{labeled Hodge-Tate weight}
$\mathbf{k} = (k_\tau)_{\tau:K \hookrightarrow E}$
is regular in the sense that each $k_\tau\in \bZ^N$ contains distinct numbers.
\end{lem}

\begin{proof}
It follows from Proposition \ref{prop:lbl-vs-colbl}.
\end{proof}

\begin{lem}
Let $K'/K$ be a finite field extension
such that $K'\subset E$.
Let $\rho:G_K \to G(E)$ be a Hodge-Tate $G$-valued representation.
We have 
$\hT(\rho|_{G_{K'}}) = \hT(\rho)$.
\end{lem}

\begin{proof}
Note that the Definition \ref{def:ht}
only makes use of $G_E$ and does not depend on~$K$.
\end{proof}

\begin{lem}
Let $\rho_1, \rho_2:G_K \to G(E)$ be two Hodge-Tate representations whose
image is abelian and consists of semisimple elements.
If $\rho_1\rho_2 = \rho_2\rho_1$, then
$\hT(\rho_1\rho_2) = \hT(\rho_1)\hT(\rho_2)$.
\end{lem}

\begin{proof}
By the previous lemma, it is harmless to shrink $G_K$ and thus we can assume $\rho_1$, $\rho_2$ both factor through a maximal torus
$T$ of $G$. 
By descent, we can assume $T$ is split.
Write $i: T\hookrightarrow G$ for the embedding of the maximal torus $T$.
Let $t_1, t_2:G_K \to T(E)$ be representations such that $i(t_1)=\rho_1$ and $i(t_2)=\rho_2$.

We have $\hT(\rho_1) = i(\hT(t_1))$ and $\hT(\rho_2) = i(\hT(t_2))$ by functoriality (Lemma \ref{lem:ht-functorial}).
So it suffices to show $\hT(t_1t_2) = \hT(t_1)\hT(t_2)$.
Since $T$ is a split torus, the general case follows from the special case $T=\Gm$.
The Hodge-Tate cocharacter of $t_1t_2: G_K \to \Gm(E)$
is completely decided by the Hodge-Tate weight of $t_1t_2$. The lemma follows because
the Hodge-Tate weight of $t_1t_2$ is the sum of the Hodge-Tate weight of $t_1$ and the
Hodge-Tate weight of $t_2$.
\end{proof}

We use the following lemma to construct Hodge-Tate regular cocharacters.

\begin{lem}
\label{lem:regular-construction}
Assume $E=K_f$ is the unramified extension of $K$ of degree $f$ inside the fixed algebraic closure $\overline{K}:=\overline{\Qp}$ of $K$.
Fix a maximal split torus $T$ of $G$.
Write $i:T \to G$ for the embedding.

For each colabel $\sigma_0: K_f \hookrightarrow \bC$, and each cocharacter
$\lambda\in X_*(G(\bC), T(\bC))$,
there exists a crystalline representation
$t:G_{K_f}\to T(K_f)$ such that
$$
\hT(i(t))^{\sigma} = 
\begin{cases}
\lambda & \text{if} \hspace{5mm} \sigma=\sigma_0, \\
\text{the trivial cocharacter} & \text{if otherwise.}
\end{cases}
$$
\end{lem}

\begin{proof}
Let $\chi_{\LT}: G_{K_f}\to \cO_{K_f}^*$ be a Lubin-Tate character.
Choose an isomorphism $T \cong \Gm^{\times r}$, $r=\rk T$.

The field $K_f$ is a subfield of $\overline{K}$ by its choice.
The composite $K_f \hookrightarrow \overline{K} \hookrightarrow \bC$
defines a canonical embedding of $K_f$ in $\bC$.
Since $K_f/K$ is a Galois extension, there exists a unique $\iota \in \Gal(K_f/K)$
such that $\sigma_0\circ \iota$ is the canonical embedding $K_f \hookrightarrow \bC_K$.

Put $t = \iota(\chi_{\LT}^{h_1}, \cdots, \chi_{\LT}^{h_r})$, $h_1,\cdots,h_r\in \bZ$.
By Lemma \ref{lem:ht-functorial} and Lemma \ref{lem:ht-lubin-tate},
$\hT(i(t))^{\sigma}$ is the trivial cocharacter if $\sigma \ne \sigma_0$.
Since the co-labeled Hodge-Tate weights of the Lubin-Tate character is $(1,0\cdots,0)$,
if we let the tuple $(h_1,\cdots,h_r)$ range over all $\bZ^{r}$,
then $\hT(i(\iota(\chi_{\LT}^{h_1}, \cdots, \chi_{\LT}^{h_r})))^{\sigma_0}$ ranges over all cocharacters in $X_*(G(\bC), T(\bC))$.
So we can choose $(h_1,\cdots,h_r)$ so that $\hT(i(\iota(\chi_{\LT}^{h_1}, \cdots, \chi_{\LT}^{h_r})))^{\sigma_0}=\lambda$.
\end{proof}

\subsection{Hodge-Tate regular lifts of quasi-semisimple mod $\varpi$ Galois representations}~

In many applications, we need Hodge-Tate regular crystalline representations.
For example, crystalline deformation rings of regular Hodge-Tate weights
have the largest dimension, which is exploited in the work \cite{EG19}.

The following lemma shows as long as a crystalline lift exists,
Hodge-Tate regular lifts also exist.

We will specialize to the case where $E = K_f$, the unramified extension of $K$ of degree $f$.

\subsubsection{Local class field theory}
Let $\Art_K: K^{\times} \to G_K^{\ab}$ be the local Artin map, which we normalize so that a uniformizer corresponds to a geometric Frobenius element.

Note that $\Art_K$ induces an isomorphism 
$$\Art_K^{-1}:\Gal(K^{\ab}/K^{\ur}) \xrightarrow{\cong} \cO_K^{\times}$$
See the paragraph after the proof of 
\cite[6.2]{Iw86} for a reference.
Denote by $r_K$ the induced map $I_K \to \cO_K^{\times}$.

\paragraph{\bf Theorem}{\cite[6.11]{Iw86}}
\label{thm:cft-ism}
Let $\sigma: K\to K'$ be an isomorphism of fields.
Then the following diagram is commutative:
$$
\xymatrix{
    K^{\times} \ar[d]^{\sigma} \ar[r]^{\Art_K} & G_K^{\ab}\ar[d]^{\sigma^*}\\
    K'^{\times} \ar[r]^{\Art_{K'}} & G_{K'}^{\ab}
}
$$
Here $\sigma^*: \tau \mapsto \sigma\tau\sigma^{-1}$.

\begin{cor}
\label{cor:rF-twist}
Let $\sigma: K\to K$ be a continuous field automorphism.
Then $r_K(\sigma \tau \sigma^{-1}) = \sigma (r_K(\tau))$ for all $\tau\in I_K$.
\end{cor}

\begin{proof}
It is an immediate consequence of Theorem \ref{thm:cft-ism}.
\end{proof}

\begin{thm}
\label{thm:ht-regular-lift}
Let $\bar\rho : G_K \to G(\bF)$ be a $G$-completely reducible representation.
Let $\kappa$ be the residue field of $K$.
Assume $\kappa\subset \bF$.

(1) There exists a Hodge-Tate regular crystalline lift
$\rho: G_K \to G(\cO_{K_f})$
for some positive integer $f$.

(2) If $G$ has a simply connected derived subgroup and $\bF$ is the splitting field of $\bar\rho|_{I_K}$ (see Theorem \ref{thm:irreducible-lift}), then $f$ can be taken as $[\bF:\kappa]$.
\end{thm}

\begin{proof}
Write $i:T\hookrightarrow G$ for the embedding of the maximal torus $T$.

We will show that as long as a crystalline lift exists, a Hodge-Tate regular crystalline also exists with the same coefficient field.
The existence of crystalline lifts is Theorem \ref{thm:irreducible-lift}.

We keep notations used in the proof of Proposition \ref{prop:lift-inertia}.
We set $\cO := \cO_{K_f}$.
Recall that $\Xi := \sum_{i=0}^{f-1}\mathpzc{w}^i\otimes \Phi^{f-1-i}_K$,
where $\Phi_K \in G_K/I_K$ is a generator of $G_K/I_K$, and $\mathpzc{w} \in W(G,T)$ is the Weyl group element which corresponds to $\bar\rho(\Phi_K)^{-1}$.
Recall that the submodule $M^0_{T,\cris} \subset M_{T,\cris}$
consists of representations $I_K \to T(\cO)$ which are extendable to $G_{K_f}$.
For each element of $u\in M^0_{T,\cris}$, choose an extension $t_u: G_{K_f} \to T(\cO)$.
The Hodge-Tate cocharacter $\hT(i(t_u))$ does not depend on the choice of $t_u$.
It makes sense to write $\hT(u)$ for $\hT(i(t_u))$ (where $t_u$ is any choice of extension).

In the proof of Proposition \ref{prop:lift-inertia}, we've shown that there exists
$v\in \Xi M^0_{T,\cris}\subset M^0_{T,\mathpzc{w},\cris}$
which is a lift of $\bar\rho|_{I_K}$.

Fix a colabel $\sigma_0: K_f \hookrightarrow \bC$.
By Lemma \ref{lem:regular-construction}, there exists a crystalline representation
$t: G_{K_f} \to T(\cO)$ such that
$\hT(i(t))^{\sigma}$ is a regular cocharacter in $X_*(G(\bC), T(\bC))$
if $\sigma=\sigma_0$, and is the trivial cocharacter if $\sigma\ne \sigma_0$.

The restriction $t|_{I_K}$ defines an element $v_0\in M^0_{T,\cris}$.
By Lemma \ref{lem:ht-functorial}, 
we have $$\hT((\mathpzc{w}\otimes 1)v_0) = \mathpzc{w}\hT(v_0)\mathpzc{w}^{-1}.$$
By Lemma \ref{lem:ht-lubin-tate} and Corollary \ref{cor:rF-twist},
we have $$\hT((1\otimes \Phi_K)v_0)^{\sigma} = \hT(v_0)^{\sigma\circ \Phi_K^{-1}}.$$
Summing up, we have
\begin{align*}
\hT(\Xi v_0)^{\sigma_0 \circ \Phi_K^{-1-i+f}} & = 
\hT(\sum_{j=0}^{f-1}\mathpzc{w}^j\otimes \Phi^{f-1-j}_K v_0)^{\sigma_0 \circ \Phi_K^{-1-i+f}} \\
&=
\prod_{j=0}^{f-1}
\hT(\mathpzc{w}^j\otimes \Phi^{f-1-j}_K v_0)^{\sigma_0 \circ \Phi_K^{-1-i+f}} \\
&=
\prod_{j=0}^{f-1}
\mathpzc{w}^j\hT(1 \otimes \Phi^{f-1-j}_K v_0)^{\sigma_0 \circ \Phi_K^{-1-i+f}}\mathpzc{w}^{-j} \\
&=
\prod_{j=0}^{f-1}
\mathpzc{w}^j\hT(v_0)^{\sigma_0 \circ \Phi_K^{-1-i+f}\circ \Phi_K^{1+j-f}}\mathpzc{w}^{-j} \\
&=
\mathpzc{w}^i\hT(v_0)^{\sigma_0}\mathpzc{w}^{-i}
\end{align*}
By Definition \ref{HT-regular}, $\Xi v_0$ is Hodge-Tate regular.

Let $C$ be a very large positive integer.
Write $N$ for the cardinality of $\bF^\times$.
Define $v' := v + C N \Xi v_0$.
Since $M_{T, \modp}$ is $N$-torsion, 
$v'$ is a lift of $\bar\rho|_{I_K}$.
We have $\hT(v') = \hT(v)\hT(\Xi v_0)^{C N}$.
Since $\hT(\Xi v_0)$ is a regular cocharacter, $\hT(v')$ is also a regular cocharacter if $C \gg 0$.

Since $\Xi M_{T, \text{cris}}^0 \subset M^0_{T,\mathpzc{w},\text{cris}}$,
we have $v + \Xi v_0 \in M_{T,\mathpzc{w},\text{cris}}$.
By Corollary \ref{cor:lift-witt},
$v'$ extends to a representation $G_K\to G(\cO)$
which is a crystalline representation lifting $\bar\rho$.
\end{proof}

\begin{appendix}
\section{Appendix: Hodge-Tate theory with coefficients}

\label{Hodge-Tate}

Let $K/\Qp$, $E/\Qp$ be finite extensions.
Assume $E$ admits an embedding of the Galois closure of $K$.
Fix an embedding $K \hookrightarrow E$.
Let $V$ be a finite dimensional $E$-vector space.
Let $\rho:G_K \to \GL(V)$ be a continuous representation.
Assume~$\rho$ is Hodge-Tate.
Let $\bC := \bC_K$ be the completed algebraic closure of $K$.
Let $\BHT := \bigoplus_{n\in\bZ}\bC(n)$ be the Hodge-Tate period ring.
Then $\BHT\otimes V := \BHT\otimes_{\Qp} V$ is a $\bC \otimes E$-module
with $G_K$-action.

Let $\sigma$ be an embedding $E \hookrightarrow \bC$.
Define 
\begin{eqnarray}
V_{\sigma} &:=&\{\sum x_i\otimes y_i\in \BHT\otimes V|
\sum \sigma(a) x_i\otimes y_i 
= \sum x_i \otimes a y_i \text{ for all }a\in E\}\nonumber\\
&=&\bigcap_{a\in E}\Ker (l_{1\otimes a} - l_{\sigma(a)\otimes 1})
\text{\hspace{9mm}(where $l_{x}$ is scalar multiplication by $x$)}
\nonumber
\end{eqnarray}
It is easy to see that

\begin{lem}
Let $L_{\sigma}\subset \bC$ be the subfield generated by $K$ and $\sigma(E)$.
\begin{itemize}
\item[(i)]
$V_{\sigma}$ is a $G_{L_{\sigma}}$-stable $\bC\otimes E$-submodule
of $\BHT\otimes V$;
\item[(ii)]
$V_{\sigma}$ is isomorphic to
$\BHT \otimes_{\sigma, E} V$
as a $G_{L_{\sigma}}$-semi-linear $\bC$-module;
\item[(iii)]
$\BHT\otimes V = \bigoplus_{\sigma: E\hookrightarrow \bC}V_{\sigma}$.
\end{itemize}
\end{lem}

Let $L$ be the Galois closure of $L_\sigma$ in $\bC$.
Write $D_{\sigma}(V) := V_{\sigma}^{G_{L}}$.
By (iii), $$
\bigoplus_{\sigma: E\hookrightarrow \bC} D_{\sigma}(V)
= (\BHT\otimes V)^{G_L}
= D_{\HT}(V)\otimes_{K}L
$$
The Hodge-Tate grading on $D_{\HT}(V)$ induces
a grading on each of $D_{\sigma}(V)$.
So $D_{\sigma}(V_\sigma)$ is a graded $L$-vector space.
We denote by $\HT^{\sigma}(V)$
 the multiset of integers $n$
 in which $n$ occurs with multiplicity
 $\dim_L \gr^n D_{\sigma}(V_\sigma)$, and call it
 the $\sigma$-\textit{co-labeled Hodge-Tate weights of $V$}.
\footnote{This is a non-standard terminology.}

\subsection{Labeled Hodge-Tate weights}

Let $\tau: K \hookrightarrow E$ be an embedding.
Define 
\begin{eqnarray}
\tilde V_{\tau} &:=&\{\sum x_i\otimes y_i\in \BHT\otimes V|
\sum a x_i\otimes y_i 
= \sum x_i \otimes \tau(a) y_i \text{ for all }a\in K\}\nonumber\\
&=&\bigcap_{a\in K}\Ker (l_{a\otimes 1} - l_{1\otimes \tau(a)})\nonumber
\end{eqnarray}

\begin{lem}
We have
$$
\tilde V_{\tau} = \bigoplus_{\sigma:E\hookrightarrow \bC, \sigma|_{\tau K} = \tau^{-1}} V_{\sigma}
$$
\end{lem}

\begin{proof}
Unravel the definitions.
\end{proof}

While $V_{\sigma}$ is only $G_{L_{\sigma}}$-stable,
$\tilde V_{\tau}$ is $G_K$-stable!
Write $\tilde D_{\tau}(V) := (\tilde V_{\tau})^{G_K}$.
We want to remind readers the usual definition
of $\tau$-labeled Hodge-Tate weights
(for example, the definition in
\cite[1.1]{GK14}).

\begin{dfn}
The multiset $\HT_{\tau}(V)$ is as follows:
an integer $n$ appears with multiplicity
$$
\dim_{E}\gr^n(D_{\HT}(V)\otimes_{E\otimes_{\Qp} K, \tau}E)
$$
\end{dfn}

\begin{lem}
We have
$
\dim_{E}\gr^n(D_{\HT}(V)\otimes_{E\otimes_{\Qp} K, \tau}E)
=\dim_E \gr^n(\tilde D_{\tau}(V))
$.
\end{lem}

\begin{proof}
It is easy to see (by unravelling the definitions) that the natural map
$$
\tilde V_{\tau}^{G_K} \hookrightarrow D_{\HT}(V)
\twoheadrightarrow D_{\HT}(V)\otimes_{E\otimes_{\Qp}K, \tau}E
$$
is injective, and $E$-linear.
So it must be an $E$-isomorphism because of the direct sum decomposition.
\end{proof}

When we divide a multiset by an integer $s$, we divide the multiplicity of all members of the multiset by $s$.
For example $\frac{1}{2}\{1,1,2,2,2,2\}=\{1,2,2\}$.

\begin{prop}
\label{prop:lbl-vs-colbl}
We have
$\HT_{\tau}(V) = \frac{1}{[E:K]}
\bigcup_{\sigma:E\hookrightarrow\bC, \sigma|_{\tau K}=\tau^{-1}}\HT^{\sigma}(V)$.
\end{prop}

\begin{proof}
Let $L$ be as before.
We have 
\begin{eqnarray}
 \tilde D_{\tau}(V)\otimes_KL &=&
\tilde V_{\tau}^{G_L} =
\bigoplus_{\sigma:E\hookrightarrow \bC, \sigma|_{\tau K} = \tau^{-1}} V_{\sigma}^{G_L} =
\bigoplus_{\sigma:E\hookrightarrow \bC, \sigma|_{\tau K} = \tau^{-1}} D_{\sigma}(V)\nonumber
\end{eqnarray}
as graded modules.
So
\begin{align*}
&\dim_E(\tilde D_{\tau}(V))
=\frac{1}{[E:K]}\dim_K(\tilde D_{\tau}(V))
=\frac{1}{[E:K]}\dim_L(\tilde D_{\tau}(V)\otimes_KL)
\\&
=\frac{1}{[E:K]}
\sum_{\sigma:E\hookrightarrow \bC, \sigma|_{\tau K} = \tau^{-1}} \dim_LD_{\sigma}(V)
\end{align*}
\end{proof}

Thus the multiset of $\tau$-labeled Hodge-Tate weights
is the average of certain multisets of $\sigma$-co-labeled Hodge-Tate weights.

\subsection{Galois twist}

The following is a convenient observation.

\begin{lem}
\label{lem:ht-res-ext}
Let $K$, $E$ be arbitrary finite extensions of $\Qp$.
Let $L/E$ be a field extension.
Let $\sigma: E \hookrightarrow \bC$ be an embedding.
Let $\tilde \sigma: L\hookrightarrow \bC$ be an embedding extending~$\sigma$.
Let $K'/K$ be a finite extension.
Then
\begin{itemize}
\item 
[(1)] $\HT^{\sigma}(\Res^{G_{K'}}_{G_K}V) = \HT^{\sigma}(V)$;
\item
[(2)]
$\HT^{\sigma}(V) = \HT^{\tilde \sigma}(V\otimes_EL)$.
\end{itemize}
Assume moreover that $E$ admits an embedding of the Galois closure of $K$.
Let $\tau:K\hookrightarrow E$ be an embedding.
Then

\begin{itemize}
\item 
[(3)]
$\HT_{\tau}(V) = \HT_{\tau}(V\otimes_EL)$.
\end{itemize}
\end{lem}

\begin{proof}
(1), (3): unravel definitions;
(2):
$\BHT\otimes_{L, \tilde \sigma}(V\otimes_EL) = (\BHT\otimes_{L, \tilde \sigma}L)\otimes_{E}V
=\BHT\otimes_{E, \sigma}V$.
\end{proof}

\begin{cor}
\label{cor:galois-twist}
Assume $E$ contains the Galois closure of $K$.
Let $\theta\in \Aut(E/\Qp)$.
Let $\tau: K\hookrightarrow E$ be an embedding.
Then

(1)
$\HT^{\sigma}(V\otimes_{E, \theta}E) = \HT^{\sigma\circ \theta}(V)$.

(2)
$\HT_{\tau}(V\otimes_{E, \theta}E) = \HT_{\theta^{-1}\circ\tau}(V)$.
\end{cor}

\begin{proof}
(1)
It is a special case of Lemma \ref{lem:ht-res-ext}(2).

(2)
By Proposition \ref{prop:lbl-vs-colbl},
\begin{eqnarray*}
\HT_{\tau}(V\otimes_{E,\theta}E) &=& 
\frac{1}{[E:K]}
\sum_{\sigma:E\hookrightarrow\bC, \sigma|_{\tau K}=\tau^{-1}}\HT^{\sigma}(V\otimes_{E,\theta}E)\nonumber\\
&=&
\frac{1}{[E:K]}
\sum_{\sigma:E\hookrightarrow\bC, \sigma|_{\tau K}=\tau^{-1}}\HT^{\sigma\circ \theta}(V)\nonumber\\
&=&
\frac{1}{[E:K]}
\sum_{\sigma:E\hookrightarrow\bC, \sigma\circ\theta^{-1}|_{\tau K}=\tau^{-1}}\HT^{\sigma}(V)\nonumber\\
&=&
\HT_{\theta^{-1}\circ\tau}(V)\nonumber\qedhere
\end{eqnarray*}
\end{proof}

\subsection{Lubin-Tate characters}

In this subsection, we want to rewrite some results of
\cite[III.A.1-III.A.5]{Se89} using the language we just developed.

\paragraph{Remark}
Proposition B.2 of \cite[Appendix B]{C11} contains a result more general than this subsection.

\subsubsection{}
Note that \textit{the cyclotomic character has Hodge-Tate weight $-1$.}

\subsubsection{Lubin-Tate characters of Galois extensions of $\Qp$}
\label{lem:LT-galois}
We start with the simpliest case.
Let $E=K/\Qp$ be a finite Galois extension.
Let $\pi$ be a uniformizer of $K$.
Let $F_\pi$ be the Lubin-Tate formal group associated to $K$ and $\pi$.
Let $\chi_K:= \chi_{K, \pi}:G_K \to \cO_E^{\times}$ be the Tate module of $F_\pi$, as is the notation of \cite{Se89}.
Then $\chi_K|_{I_K} = r_K^{\otimes -1}$ (see subsection \ref{cor:rF-twist}).
(So $r_K$ is crystalline.)

\paragraph{\bf Lemma}

Let $\sigma_1\in\Gal(K/\Qp)$.
Then a $\sigma$-co-labeled Hodge-Tate weight
of $\sigma_1\circ \chi_K$ is $-1$ if $\sigma=\sigma_1^{-1}$, and $0$ if otherwise.

\begin{proof}
See \cite[Thm 2, III.A.5]{Se89}
and \cite[Prop III.A.4]{Se89}.
Note that
\begin{itemize}
\item Serre's $K$ and $E$ are reversed,
\item Galois hypothesis is required by \cite[III.A.3(b)]{Se89},
\item Serre's $W_{\sigma}$ is our $\gr^0V_{\sigma}$.
\end{itemize}
\end{proof}

\begin{lem}
\label{lem:HT-LT}
Now suppose $E=K/\Qp$ is not necessarily Galois.
A $\sigma$-co-labeled Hodge-Tate weight of
$\chi_K$ is $-1$ if $\sigma=\Id$ 
\footnote{More precisely the tautological embedding of $E$ in $\bC$}, and $0$ if otherwise.
\end{lem}

\begin{proof}
Choose a Galois closure $L$ of $K$ over $\Qp$. Consider
$$
\xymatrix{
    G_L^{\ab} \ar[d]^{\operatorname{rec}} \ar[r] &
    G_K^{\ab} \ar[d]^{\operatorname{rec}} \\
    L^{\times} \ar[r]^{N_{L/K}} & K^{\times}
}
$$
By local class field theory,
$\chi_K|_{G_L}=N_{L/K}\circ \chi_L
=\prod_{\sigma\in \Gal(L/K)}\sigma\circ \chi_L$.
By Lemma \ref{lem:LT-galois},
for $\tau\in \Gal(L/\Qp)$,
$$
\HT^{\tau}(\chi_K|_{G_L})=
\begin{cases}
-1 & \text{if $\tau$ fixes $K$} \\
0 & \text{if otherwise}
\end{cases}
$$
Now apply Lemma \ref{lem:ht-res-ext}(1), (2) to conclude.
\end{proof}

\begin{lem}
\label{lem:ht-lubin-tate}
Let $K/\Qp$ be a finite extension, and let $E/\Qp$ be a finite extension admitting $\iota: K \hookrightarrow E$.

(1)
For each $\sigma:E\hookrightarrow \bC$,
the $\sigma$-co-labeled Hodge-Tate weight
of $\iota\circ \chi_K$ is $-1$ if
$\sigma\circ \iota = \Id_K$, and $0$ if otherwise.

(2)
Suppose further $E$ admits an embedding
of the normal closure of $K$.
Then for each $\sigma:K\to E$,
the $\sigma$-labeled Hodge-Tate weight of $\iota\circ \chi_K$ is $-1$ if
$\sigma=\iota$, and $0$ if otherwise.
\end{lem}

\begin{proof}
(1) We have
\begin{align*}
\HT^{\sigma}(\iota\circ\chi_K) &
=\HT^{\sigma\circ \iota}(\chi_K) & \text{By Lemma \ref{lem:ht-res-ext}(2)} \\
&=
\begin{cases}
 -1 & \text{if $\sigma\circ\iota=\Id_K$}\\
 0 & \text{if otherwise}
\end{cases}
&\text{By Lemma \ref{lem:HT-LT}}
\end{align*}

(2)
Follows from Proposition \ref{prop:lbl-vs-colbl} and (1).
\end{proof}

\begin{lem}
\label{lem:ind-calc}
Let $K/\Qp$ be a finite extension.
Let $L/K$ be an unramified extension in $\bC$.
Let $L'$ be the Galois closure of $L$
over $\Qp$.
Let $\iota:K\hookrightarrow L'$ be the tautological embedding.
Let $\Phi_K\in G_K$ be a lift of a topological generator of $G_K/I_K$.
Let $d = [L:K]$.
Then 

(1)
Let $\sigma: L \hookrightarrow \bC$. Then
$$
\HT^{\sigma}(\Ind_{G_L}^{G_K} (\chi_L)) = 
\HT^{\sigma}(\chi_L) \cup \HT^{\sigma}(\Phi_K\circ \chi_L) \cup \dots \cup \HT^{\sigma}(\Phi_K^{d-1}\circ\chi_L)
$$

(2)
Let $\tau: K\hookrightarrow \bC$. Then
$$
\HT_{\tau}(\Ind_{G_L}^{G_K}(\iota\circ \chi_L)) = 
\begin{cases}
\{0, \dots, 0, -1\}& \text{if }\tau\text{ is the canonical embedding, and}\\
\{0, \dots, 0, 0\} &\text{if otherwise.}
\end{cases}
$$
\end{lem}

\begin{proof}
(1)
Follows from Lemma \ref{lem:ht-res-ext} and Corollary \ref{cor:rF-twist}.

(2) We have
\begin{align*}
&\HT_{\tau}(\Ind_{G_L}^{G_K}(\iota\circ \chi_L)) =
\frac{1}{[L':K]}\bigcup_{\tilde \sigma: L'\hookrightarrow \bC, \tilde \sigma|_{\tau K}\circ \tau = id}
\HT^{\tilde \sigma}(\Ind_{G_L}^{G_K}(\iota\circ \chi_L))\nonumber\\
 = &
\frac{1}{[L':K]}\bigcup_{\tilde \sigma: L'\hookrightarrow \bC, \tilde \sigma|_{\tau K}\circ \tau = id,
\sigma := \tilde \sigma|_L}
\HT^{\sigma}(\Ind_{G_L}^{G_K}( \chi_L))\nonumber\\
 = &
\frac{1}{[L:K]}\bigcup_{\sigma: L\hookrightarrow \bC, \sigma|_{\tau K}\circ \tau = id}
\HT^{ \sigma}(\Ind_{G_L}^{G_K}(\chi_L))\nonumber\\
 = &
\frac{1}{[L:K]}\bigcup_{\sigma: L\hookrightarrow \bC, \sigma|_{\tau K}\circ \tau = id}
\HT^{\sigma}(\chi_L) \cup \HT^{\sigma}(\Phi_K\circ \chi_L) \cup \dots \cup \HT^{\sigma}(\Phi_K^{d-1}\circ\chi_L)\nonumber\\
 = &
\frac{1}{[L:K]}\bigcup_{\sigma: L\hookrightarrow \bC, \sigma|_{\tau K}\circ \tau = id}
\bigcup_{k = 0}^{d-1}
\delta_{\sigma, \iota\circ \Phi_K^k}\nonumber
\end{align*}
Here $\delta_{X, Y}$ is $\{-1\}$ if $X=Y$ and is $\{0\}$ if otherwise.
Since $\Phi^k_K$ is the identity when restricted on $K$,
the last line is $0$ unless $\tau$ is the canonical embedding;
in this case, the last line becomes
$
\frac{1}{[L:K]}\bigcup_{j = 0}^{d-1}
\bigcup_{k = 0}^{d-1}
\delta_{\iota\circ\Phi_K^j, \iota\circ \Phi_K^k} = \{0, \dots, 0, -1\}
$.
\end{proof}

\end{appendix}


\end{document}